\documentclass[12pt]{article}%
\usepackage{amsmath}
\usepackage{amsfonts}
\usepackage{amssymb}
\usepackage{graphicx}%
\providecommand{\U}[1]{\protect\rule{.1in}{.1in}}
\newtheorem{theorem}{Theorem}

\begin{document}

\title{The class three groups of order $p^{9}$ with exponent $p$}
\author{Michael Vaughan-Lee}
\date{February 2017}
\maketitle

\begin{abstract}
In this note we give a complete list of all class two groups of prime exponent
with order $p^{k}$ for $k\leq8$. For every group in this list we are able to
show that the number of conjugacy classes of the group is polynomial in $p$,
and that the order of the automorphism group is also polynomial in $p$. In
addition, for each group in the list, the number of class 3 descendants of
order $p^{9}$ with exponent $p$ is PORC. So the total number of class 3 groups
of order $p^{9}$ with exponent $p$ is PORC.

\end{abstract}

\section{Introduction}

In 1960 Graham Higman wrote two immensely important and significant papers
\cite{higman60}, \cite{higman60b}. In these papers he conjectured that if we
set $f(p,n)$ equal to the number of groups of order $p^{n}$, then $f(p,n)$ is
a PORC function of $p$ for each positive integer $n$ --- this is known has
Higman's PORC conjecture. (A function $f(p)$ is said to be PORC,
\textbf{P}olynomial \textbf{O}n \textbf{R}esidue \textbf{C}lasses, if there is
a finite set of polynomials in $p$, $g_{1}(p),g_{2}(p),\ldots,g_{k}(p)$, and a
fixed integer $M$, such that for each prime $p$, $f(p)=g_{i}(p)$ for some $i$
($1\leq i\leq k$), with the choice of $i$ depending on the residue class of
$p\operatorname{mod}M$.) He also proved that the number of $p$-class two
groups of order $p^{n}$ is a PORC function of $p$ for each $n$. At the time
Higman made his conjecture it was known to hold true for $n\leq5$, and it has
since been proved correct for $n\leq7$ (see \cite{newobvl}, \cite{obrienvl2}).
For example, {for }$p\geq5$ the number of groups of order $p^{6}$ is{%
\[
3p^{2}+39p+344+24\gcd(p-1,3)+11\gcd(p-1,4)+2\gcd(p-1,5).
\]
The number of groups of order }$p^{8}$ with exponent $p$ has also been shown
to be PORC (see \cite{vlee2014}). For reasons that I will explain below I
became interested in the question of whether or not the number of class 3
groups of order $p^{9}$ with exponent $p$ is PORC. These groups are all
immediate descendants of class two groups of exponent $p$ with order dividing
$p^{8}$. (A $p$-group $H$ is said to be an immediate descendant of a $p$-group
$G$ if $G\cong H/N$, where $N$ is the last non-trivial term of the lower
$p$-central series of $H$. See \cite{newobvl} for details.) In \cite{vlee2015}
I gave a complete list of the 70 class two groups of exponent $p$ ($p>2$),
together with a description of their automorphism groups. In this paper I
consider each of these 70 groups in turn, and for each group give the number
of class 3 descendants of order $p^9$ with exponent $p$. (The numbers 
given are valid for all $p \ge 5$.) In each case the number of descendants
is PORC, and so:

\begin{theorem}
The number of class three groups of order $p^{9}$ with exponent $p$ is PORC.
\end{theorem}

However there are finite $p$-groups for which the number of immediate
descendants of given order is \emph{not} PORC. For example du Sautoy and
Vaughan-Lee \cite{dusautoyvl} give an example of a class two group $G$ with
exponent $p$ and order $p^{9}$ with the following properties:

\begin{enumerate}
\item The number of conjugacy classes of $G$ is not PORC,

\item The order of the automorphism group of $G$ is not PORC,

\item The number of immediate descendants of $G$ with exponent $p$ and order
$p^{10}$ is not PORC.
\end{enumerate}

This example shows that at least one group of order $p^{9}$ has a non-PORC
number of descendants of order $p^{10}$. However this does not disprove
Higman's PORC conjecture for groups of order $p^{10}$. If we were to follow
the methods of \cite{newobvl} and \cite{obrienvl2} in an attempt to settle the
conjecture for $p^{10}$ we would need to draw up a complete and irredundant
list of the groups with order dividing $p^{9}$, and compute the number of
descendants of order $p^{10}$ of each group in the list. There are likely to
be other groups of order $p^{9}$ with a non-PORC number of descendants of
order $p^{10}$, so the grand total of groups of order $p^{10}$ could be PORC,
with individual non-PORC summands cancelling out. Nevertheless du Sautoy and I
confidently asserted: \textquotedblleft\emph{It seems likely that there are
other groups of order }$p^{9}$\emph{\ with a non-PORC number of immediate
descendants of order }$p^{10}$\emph{, and so it is possible that the grand
total is PORC, even though not all of the summands are PORC. The authors' own
view is that this is extremely unlikely.}\textquotedblright

It is still an open question whether the number of groups of order $p^{10}$ is
PORC, but nevertheless this assertion is now looking rather foolish! Seungjai
Lee \cite{seungjai} has recently found a class two group $H$ with exponent $p$
and order $p^{8}$ with the following properties:

\begin{enumerate}
\item The number of conjugacy classes of $H$ is not PORC,

\item The order of the automorphism group of $H$ is not PORC,

\item The number of immediate descendants of $H$ with exponent $p$ and order
$p^{9}$ is not PORC.
\end{enumerate}

However, in this case we now know that the total number of class three groups
of order $p^{9}$ with exponent $p$ is PORC. In fact, in apparent contradiction
to Lee's result, for every single group $G$ in my list of 70 class two groups
of exponent $p$ with order dividing $p^{8}$:

\begin{enumerate}
\item The number of conjugacy classes of $G$ is a polynomial in $p$,

\item The order of the automorphism group of $G$ is a polynomial in $p$,

\item The number of class 3 descendants of $G$ with order $p^{9}$ and exponent
$p$ is PORC.
\end{enumerate}

The explanation for this \textquotedblleft contradiction\textquotedblright\ is
that Lee's group $H$ is really a family of groups, one for each $p$, as are
the groups in my list. For any given $p$, Lee's group $H$ must lie in one of
my families, but it does not have to lie in the same family for every $p$. It
turns out that there are four (families of) groups in my list $A$, $B $, $C$,
$D$. If $p=3$ then $H\cong A$, if $p=2\operatorname{mod}3$ then $H\cong B$, if
$p=1\operatorname{mod}3$ and $t^{3}-2$ has no roots in GF$(p)$ then $H\cong
C$, and if $p=1\operatorname{mod}3$ and $t^{3}-2$ has three roots in GF$(p)$
then $H\cong D$. The non-PORC properties of Lee's group arise from the fact
that the number of roots of $t^{3}-2$ over GF$(p)$ is not PORC. (There is a
deep theorem in Number Theory that asserts that if we have a polynomial with
integer coefficients then the number of roots of the polynomial over GF$(p)$
is a PORC function of $p$ if and only if the Galois group of the polynomial
over the rationals is abelian. The Galois group of $t^{3}-2$ over the
rationals is the symmetric group of degree 3.)

\section{Methods}

The groups listed here correspond to Lie algebras over GF$(p)$ under the
Lazard correspondence. All the presentations are for class two groups of
exponent $p$ and the presentations can also be read as presentations for class
two Lie algebras over GF$(p)$, except that products have to read as sums. For
example group 8.4.4 below has presentation%
\[
\langle a,b,c,d\,|\,[d,b][c,a],\,[d,c][b,a]^{\omega}\rangle.
\]
This corresponds to the Lie algebra with presentation%
\[
\langle a,b,c,d\,|\,[d,b]+[c,a],\,[d,c]+\omega\lbrack b,a]\rangle
\]
using $[,]$ to denote the Lie product, or (as I prefer)%
\[
\langle a,b,c,d\,|\,db+ca,\,dc+\omega ba\rangle
\]
using juxtaposition to denote the Lie product. (Here, and throughout this
paper, $\omega$ denotes a fixed integer which is primitive modulo $p$.) All my
calculations are carried out in Lie algebras, and much of the detailed
information in this paper is given in Lie algebra notation. I make no apology
for mixing up group and Lie algebra notation in this way. Theorem 1 is both a
theorem in groups and in Lie algebras. A group theorist who wants to look at
the groups in this paper would probably prefer to see presentations in group
notation, and so this is the way I have given them. But the calculations are
all in Lie algebras, and it would be awkward to translate their details back
into group notation.

I have computed the number of descendants of each group using the Lie algebra
version of the $p$-group generation algorithm (see Newman \cite{Newman77} and
O'Brien \cite{OBrien90}). Let $L$ be a class two Lie algebra over GF$(p)$. The
first step is to compute the $p$-covering algebra which is the largest class 3
Lie algebra $K$ with an ideal $M$ such that

\begin{enumerate}
\item $K/M\cong L$,

\item $M\leq\zeta(K)\cap K^{2}$.
\end{enumerate}

(Here $\zeta(K)$ denotes the centre of $K$.) It is quite easy to construct $K$
with the prime $p$ symbolic. (So in effect we construct a family of covering
algebras, one for each $p$.) The ideal $M$ is called the $p$-multiplier of
$L$, and $N=K^{3} \le M$ is called the nucleus of $L$. An allowable subalgebra
of $M$ is a proper subalgebra $S$ such that $S+N=M$. The immediate descendants
of $L$ are the Lie algebras $L/S$ where $S$ is an allowable subalgebra. The
automorphism group of $L$ acts on $M$, and two descendants $L/S$, $L/T$ are
isomorphic if and only if $S$ and $T$ are in the same orbit under the action
of the automorphism group. Since we are only calculating descendants of
dimension 9, if $L$ has dimension $n$ then we can restrict our attention to
allowable subalgebras of $M$ of codimension $9-n$. To compute the number of
immediate descendants of $L$ of dimension 9 we need to compute the number of
orbits of allowable subspaces of $M$ of codimension $9-n$ under the action of
the automorphism group. We mostly do this in two steps. If $S$ is an allowable
subalgebra of codimension $9-n$ in $M$, then $S\cap N$ has codimension $9-n$
in $N$. So first we compute a set of representatives for the orbits of
subspaces of codimension $9-n$ in $N$. We can restrict our attention to
allowable subalgebras of $M$ such that $S\cap N$ is one of these
representatives. If $T$ is one of these representatives, and if $S\cap N=T$,
then%
\[
K/S\cong(K/T)/(S/T).
\]
Furthermore the orbits of allowable subalgebras $S$ with $S\cap N=T$
correspond to orbits under the stabilizer of $T$ of subalgebras of codimension
$9-n$ in $M/T$ which have trivial intersection with $N/T$. For most of
the Lie algebras with descendants of dimension 9 we give the dimensions of $M$
and $N$, and a set of representatives for the orbits of subspaces of $N$ of
codimension $9-n$. Then for each such representative $T$ we give the number of
descendants $K/S$ where $S\cap N=T$.

In a few cases I use Higman's ideas to compute the number of orbits of
subalgebras of $M$ of codimension $9-n$, and then obtain the number of
descendants by counting the number of orbits of non-allowable subalgebras.

\section{Order $p^{3}$}

\noindent\textbf{Group 3.2.1}%

\[
\langle a,b\rangle
\]

(There are no relations here, since this is the free class two group of 
exponent $p$ on two generators. As mentioned above, all the presentations are
for class two groups of exponent $p$.) The number of conjugacy classes is
$p^{2}+p-1$, and the automorphism group has order $(p^{2}-1)(p^{2}-p)p^{2}$.
There are no class 3 descendants of order $p^{9}$.

\section{Order $p^{4}$}

\noindent\textbf{Group 4.3.1}%

\[
\langle a,b,c\,|\,[c,a],\,[c,b]\rangle
\]

The number of conjugacy classes is $p^{3}+p^{2}-p$, and the automorphism group
has order $(p-1)(p^{2}-1)(p^{2}-p)p^{5}$. There are no class 3 descendants of
order $p^{9}$.

\section{Order $p^{5}$}

\setcounter{subsection}{2}

\subsection{Three generator groups}

\noindent\textbf{Group 5.3.1}%

\[
\langle a,b,c\,|\,[c,b]\rangle
\]

The number of conjugacy classes is $2p^{3}-p$, and the automorphism group has
order $(p-1)(p^{2}-1)(p^{2}-p)p^{8}$. This group has 12 class 3 descendants of
order $p^{9}$ and exponent $p$.

\subsection{Four generator groups}

\noindent\textbf{Group 5.4.1}%

\[
\langle a,b,c,d\,|\,[c,a],\,[c,b],\,[d,a],\,[d,b],\,[d,c]\rangle
\]

The number of conjugacy classes is $p^{4}+p^{3}-p^{2}$, and the automorphism
group has order $(p^{2}-1)^{2}(p^{2}-p)^{2}p^{8}$. There are no class 3
descendants of order $p^{9}$.

\bigskip\noindent\textbf{Group 5.4.2}%
\[
\langle a,b\rangle\times_{\lbrack b,a]=[d,c]}\langle c,d\rangle
\]

The number of conjugacy classes is $p^{4}+p-1$, and the automorphism group has
order $(p^{5}-p)(p^{5}-p^{4})(p^{3}-p)p^{2}$. There are no class 3 descendants
of order $p^{9}$.

\section{Order $p^{6}$}

\setcounter{subsection}{2}

\subsection{Three generator groups}

\noindent\textbf{Group 6.3.1}%

\[
\langle a,b,c\rangle
\]

The number of conjugacy classes is $p^{4}+p^{3}-p$, and the automorphism group
has order $(p^{3}-1)(p^{3}-p)(p^{3}-p^{2})p^{9}$. This group has%
\begin{align*}
&  p^{7}+p^{6}+3p^{5}+5p^{4}+11p^{3}+(22+2\gcd(p-1,3))p^{2}\\
&  +(46+5\gcd(p-1,3)+\gcd(p-1,4))p\\
&  +60+6\gcd(p-1,3)+2\gcd(p-1,4)
\end{align*}
class 3 descendants of order $p^{9}$ and exponent $p$. In this group the
$p$-multiplier $M$ is equal to the nucleus $N$, and so the number of
descendants is equal to the number of orbits of subspaces of $N$ of
codimension 3.

\subsection{Four generator groups}

\noindent\textbf{Group 6.4.1}%

\[
\langle a,b,c,d\,|\,[c,b],\,[d,a],\,[d,b],\,[d,c]\rangle
\]

The number of conjugacy classes is $2p^{4}-p^{2}$ and the order of the
automorphism group is $(p-1)^{2}(p^{2}-1)(p^{2}-p)p^{13}$. The nucleus has
dimension 5. We can assume that one of the following 16 sets of relations
holds:%
\begin{align*}
bab &  =bac=0,\\
bac &  =cac+bab=0,\\
bac &  =cac+\omega bab=0,\\
caa &  =bab=0,\\
caa &  =bab-cac=0,\\
caa &  =\text{ }bab-\omega cac=0,\\
caa &  =bac=0,\\
caa &  =cac=0,\\
bac &  =caa-bab=0,\\
cac &  =caa-bab=0,\\
baa &  =caa=0,\\
baa &  =caa-bab=0,\\
baa &  =caa-cac=0,\\
baa &  =caa-bab+cac=0,\\
baa &  =caa-bab+\omega cac=0,\\
baa-cac &  =caa-\alpha bab-cac=0,
\end{align*}
for some $\alpha$ with $x^{3}+\alpha^{-1}x-\alpha^{-1}$ irreducible over
GF$(p)$. (Throughout this paper $\omega$ denotes a fixed integer which is
primitive modulo $p$.) The number of immediate descendants of order $p^{9}$
with exponent $p$ in the 16 cases are as follows:%
\begin{align*}
&  7p+50+\gcd(p-1,4),\\
&  \frac{1}{2}(3p^{2}+17p+80+\gcd(p-1,4)),\\
&  \frac{1}{2}(3p^{2}+11p+30+\gcd(p-1,4)),\\
&  p^{4}+3p^{3}+7p^{2}+20p+69+2\gcd(p-1,3)+\gcd(p-1,4),\\
&  \frac{1}{2}(p^{5}+2p^{4}+4p^{3}+11p^{2}+23p+31)+\gcd(p-1,4),\\
&  \frac{1}{2}(p^{5}+2p^{4}+4p^{3}+11p^{2}+23p+31)+\gcd(p-1,4),\\
&  p^{3}+3p^{2}+12p+55+6\gcd(p-1,3)+3\gcd(p-1,4)+\gcd(p-1,5),\\
&  2p^{2}+9p+35+\gcd(p-1,3),\\
&  p^{4}+2p^{3}+5p^{2}+10p+17+(p+5)\gcd(p-1,3)\\
&  \;\;+(p+3)\gcd(p-1,4)+\gcd(p-1,5),\\
&  p^{3}+2p^{2}+8p+17+(p+3)\gcd(p-1,3)+\gcd(p-1,4),\\
&  p^{2}+3p+23+\gcd(p-1,3),\\
&  p^{4}+2p^{3}+4p^{2}+8p+9+(p+2)\gcd(p-1,3)+\gcd(p-1,4),\\
&  p^{5}+2p^{4}+3p^{3}+5p^{2}+10p+9+(p+1)\gcd(p-1,3),\\
&  \frac{1}{6}(p^{6}+p^{5}+p^{4}+4p^{3}+9p^{2}+7p+13+2(p+1)\gcd(p-1,3)),\\
&  \frac{1}{2}(p^{6}+p^{5}+p^{4}+2p^{3}+3p^{2}+3p+5),\\
&  \frac{1}{3}(p^{6}+p^{5}+p^{4}+p^{3}+3p^{2}+p+4+2(p+1)\gcd(p-1,3)).
\end{align*}

\bigskip\noindent\textbf{Group 6.4.2}%

\[
\langle a,b,c,d\,|\,[c,b],\,[d,a],\,[d,b]=[b,a],\,[d,c]\rangle
\]

The number of conjugacy classes is $p^{4}+2p^{3}-p^{2}-2p+1$ and the order of
the automorphism group is $2(p^{2}-1)^{2}(p^{2}-p)^{2}p^{8}$. The nucleus has
dimension 4, and is spanned by $baa$, $bab$, $caa$, $cac$. There are two
orbits of the automorphism group on the one dimensional subspaces of the
nucleus, and we take the (single) relation on the nucleus to be $cac=0$ or
$cac=bab$. The number of class 3 descendants of order $p^{9}$ with exponent $p
$ in the two cases is%
\begin{align*}
&  (p+11)\gcd(p-1,3)+26,\\
&  \frac{1}{2}(p^{3}+3p^{2}+11p+26)\gcd(p-1,3)+\frac{1}{2}(p^{2}+3p+29)\\
&  +\frac{1}{2}\gcd(p-1,4)+\gcd(p-1,9)+\gcd(p-1,12).
\end{align*}

\bigskip\noindent\textbf{Group 6.4.3}%

\[
\langle a,b,c,d\,|\,[c,b],\,[d,a],\,[d,b]=[c,a],\,[d,c]\rangle
\]

The number of conjugacy classes is $p^{4}+p^{3}-p$ and the order of the
automorphism group is $(p-1)(p^{2}-1)(p^{2}-p)p^{12}$. The nucleus has
dimension 4, and is spanned by $baa$, $bab$, $bac$ and $bad$. There are two
orbits of the automorphism group on the one dimensional subspaces of the
nucleus, and we take the (single) relation on the nucleus to be $bab=0$ or
$bad=0$. The number of class 3 descendants of order $p^{9}$ with exponent $p$
in the the first case is%
\[
p^{2}+6p+32+5\gcd(p-1,3)+2\gcd(p-1,4),
\]
and in the second case it is%
\begin{gather*}
75\text{ if }p=5,\\
p^{2}+5p+22+\gcd(p-1,3)+\gcd(p-1,5)\text{ if }p>5.
\end{gather*}

\bigskip\noindent\textbf{Group 6.4.4}%

\[
\langle a,b,c,d\,|\,[c,b],\,[d,a],\,[d,b]=[c,a],\,[d,c]=[b,a]^{\omega}\rangle
\]

The number of conjugacy classes is $p^{4}+p^{2}-1$ and the order of the
automorphism group is $2(p^{4}-1)(p^{4}-p^{2})p^{8}$. The nucleus has
dimension 4, and is spanned by $baa$, $bab$, $bac$ and $bad$. But there is
only one orbit of the automorphism group on the the one dimensional subspaces
of the nucleus, and so we can take the single relation on the nucleus to be
$bad=0$. The number of class 3 descendants of order $p^{9}$ with exponent $p$
is%
\[
\frac{1}{2}(4p^{3}+5p^{2}+23p+35)-\frac{1}{2}(p^{3}+p^{2}+5p+4)\gcd
(p-1,3)+\frac{1}{2}\gcd(p-1,4).
\]

\subsection{Five generator groups}

\noindent\textbf{Group 6.5.1}%

\[
\langle a,b\rangle\times\langle c\rangle\times\langle d\rangle\times\langle
e\rangle
\]

The number of conjugacy classes is $p^{5}+p^{4}-p^{3}$ and the order of the
automorphism group is $(p^{2}-1)(p^{2}-p)(p^{3}-1)(p^{3}-p)(p^{3}-p^{2}%
)p^{11}$. There are no class 3 descendants of order $p^{9}$.

\bigskip\noindent\textbf{Group 6.5.2}%

\[
\langle a,b\rangle\times_{\lbrack b,a]=[d,c]}\langle c,d\rangle\times\langle
e\rangle
\]

The number of conjugacy classes is $p^{5}+p^{2}-p$ and the order of the
automorphism group is $(p^{6}-p^{2})(p^{6}-p^{5})(p^{4}-p^{2})p^{3}(p^{2}-p)$.
There are no class 3 descendants of order $p^{9}$.

\section{Order $p^{7}$}

\setcounter{subsection}{3}

\subsection{Four generator groups}

\noindent\textbf{Group 7.4.1}%

\[
\langle a,b,c,d\,|\,[c,b],\,[d,b],\,[d,c]\rangle
\]

The number of conjugacy classes is $p^{5}+p^{4}-p^{2}$, and the automorphism
group has order $(p-1)(p^{3}-1)(p^{3}-p)(p^{3}-p^{2})p^{15}$. The nucleus, $N
$, has dimension 9, and is spanned by%
\[
bab,bac,bad,cac,cad,dad,baa,caa,daa.
\]
The subspace%
\[
B=\langle bab,bac,bad,cac,cad,dad\rangle
\]
of the nucleus is invariant under the automorphism group. The nucleus has%
\[
p^{2}+4p+76+3\gcd(p-1,3)+2\gcd(p-1,4)
\]
orbits of subspaces of dimension 7, and we organize these according to their
intersection with $B$. This intersection will have dimension 4, 5 or 6.

There is one orbit of subspaces of $N$ of dimension 7 containing $B$, and an
orbit representative is%
\[
\langle bab,bac,bad,cac,cad,dad,baa\rangle.
\]
There are $p+8$ descendants of group 7.4.1 of order $p^{9}$ with exponent $p$
satisfying relations%
\[
bab=bac=bad=cac=cad=dad=baa=0.
\]

There are 4 orbits of subspaces of dimension 5 in $B$. Orbit representatives
are%
\begin{align*}
&  \langle bab,bac,bad,cac,cad\rangle,\\
&  \langle bab,bac,bad,cac,dad\rangle,\\
&  \langle bab,bac,bad,\omega cac+dad,cad\rangle,\\
&  \langle bab,bac,cac+bad,cad,dad\rangle.
\end{align*}
The corresponding orbits are (respectively) the intersections with $B$ of 4,
5, 3 and 3 orbits of subspaces of dimension 7 in $N$. The 4 orbit
representatives corresponding to the first of these four are%
\begin{align*}
&  \langle bab,bac,bad,cac,cad,baa,caa\rangle,\\
&  \langle bab,bac,bad,cac,cad,baa,caa-dad\rangle,\\
&  \langle bab,bac,bad,cac,cad,baa,daa\rangle,\\
&  \langle bab,bac,bad,cac,cad,baa-dad,daa\rangle,
\end{align*}
and these yield (respectively) $5$, $p+5$, $10$, $8$ descendants of order
$p^{9}$ with exponent $p$. For the other 3 orbit representatives of subspaces
of dimension 5 in $B$ we obtain the following orbit representatives for
subspaces dimension 7 in $B$:%
\begin{align*}
&  \langle bab,bac,bad,cac,dad,baa,caa\rangle,\\
&  \langle bab,bac,bad,cac,dad,baa-cad,caa\rangle,\\
&  \langle bab,bac,bad,cac,dad,baa,caa+daa\rangle,\\
&  \langle bab,bac,bad,cac,dad,baa-cad,caa+daa\rangle,\\
&  \langle bab,bac,bad,cac,dad,caa,daa\rangle,
\end{align*}
yielding 5, 5, $(p+7)/2$, $(p+7)/2$, 6 descendants respectively;%

\begin{align*}
&  \langle bab,bac,bad,\omega cac+dad,cad,baa,caa\rangle,\\
&  \langle bab,bac,bad,\omega cac+dad,cad,baa-cac,baa\rangle,\\
&  \langle bab,bac,bad,\omega cac+dad,cad,caa,daa\rangle,
\end{align*}
yielding $(p+7)/2$, $(p+7)/2$, 4 descendants respectively, and:%
\begin{align*}
&  \langle bab,bac,cac+bad,cad,dad,baa,caa\rangle,\\
&  \langle bab,bac,cac+bad,cad,dad,baa,daa\rangle,\\
&  \langle bab,bac,cac+bad,cad,dad,caa,2\omega baa+daa\rangle,
\end{align*}
yielding $6$, $p+5$, $p+3$ descendants respectively.

The 15 orbits of subspaces of $B$ of dimension 4 have representatives%
\begin{align*}
&  \langle bab,bac,bad,cac\rangle,\\
&  \langle bab,bac,bad,cad\rangle,\\
&  \langle bab,bac,bad,\alpha cac+cad-\alpha dad\rangle \; (1+4\alpha^2 \mathrm{\; not\; square}),\\
&  \langle bab,bac,cac,dad\rangle,\\
&  \langle bab,bac,cad,dad\rangle,\\
&  \langle bab,bac,cac-dad,cad\rangle,\\
&  \langle bab,bac,\omega cac-dad,cad\rangle,\\
&  \langle bab,bac,bad-cac,cad\rangle,\\
&  \langle bab,bac,bad-cac,dad\rangle,\\
&  \langle bab,cac,cad-bad,dad-bad\rangle,\\
&  \langle bab,bad-cac,cad,\omega bad+dad\rangle,\\
&  \langle bab,bad-cac,cad,\beta bac+bad-dad\rangle,\\
&  \langle bac,bad,cad,\omega cac+dad\rangle,\\
&  \langle bac,bad,cad,bab+\omega cac-dad\rangle,\\
&  \langle bac,bad,bab-cac+\gamma cad,bab-\gamma cad-dad\rangle,
\end{align*}
where $\beta$ is chosen so that $x^{3}+x-\beta$ is irreducible over GF$(p)$,
and where $\gamma$ is chosen so that $1+\gamma^{2}$ is not square modulo $p$.

The first of these 15 subspaces gives 6 orbits of subspaces of dimension 7 in
$N$, with representatives of the form%
\[
\langle bab,bac,bad,cac,baa-xcad-ydad,caa-zcad-tdad\rangle.
\]
These 6 representatives together yield 20 descendants of order $p^{9}$ with
exponent $p$.

The second of these 15 subspaces gives $7+(\gcd(p-1,3)-1)/2$ orbits of
subspaces of dimension 7 in $N$, with representatives of the form%
\[
\langle bab,bac,bad,cad,baa-xcac-ydad,caa-zdad,daa-tcac\rangle.
\]
Between them these representatives yield $\frac{1}{2}(p+35+(p+9)\gcd(p-1,3)) $
descendants of order $p^{9}$ with exponent $p$.

The third gives $5-(\gcd(p-1,3)-1)/2$ orbits of subspaces of dimension 7 in
$N$, with representatives of the form%
\[
\langle bab,bac,bad,\omega cac+cad-\omega
dad,baa-xcac-ydad,caa-zcac-tdad,daa\rangle.
\]
Between them these representatives yield $\frac{1}{2}(5p+25-(p+3)\gcd(p-1,3))$
descendants of order $p^{9}$ with exponent $p$.

I haven't forgotten the fourth!!!

The fifth gives $5+\gcd(p-1,3)$ orbits of subspaces of dimension 7 in $N$,
with representatives of the form%
\[
\langle bab,bac,cad,dad,baa-xcac,caa-ybad,daa-zcac\rangle.
\]
Each of these representatives yields exactly 2 descendants of order $p^{9}$
with exponent $p$.

The eighth gives 4 orbits of subspaces of dimension 7 in $N$, with
representatives of the form%
\[
\langle bab,bac,bad-cac,cad,baa-xdad,caa-ydad,daa-zdad\rangle.
\]
Each of these representatives yields exactly 2 descendants of order $p^{9}$
with exponent $p$.

The thirteenth gives $3+\gcd(p-1,3)$ orbits of subspaces of dimension 7 in
$N$, with representatives of the form%
\[
\langle bac,bad,cad,\omega cac+dad,baa-xcac,caa-ybab,daa-zbab\rangle.
\]
Each of these representatives yields exactly 2 descendants of order $p^{9}$
with exponent $p$.

So far we have accounted for $46+2\gcd(p-1,3)$ of the%
\[
p^{2}+4p+76+3\gcd(p-1,3)+2\gcd(p-1,4)
\]
orbits of subspaces of dimension 7 in the nucleus $N$. Each of the remaining
orbits contribute one descendant of order $p^{9}$ with exponent $p$. So
the total number of descendants is
\[
p^2+13p+188+8\gcd(p-1,3)+2\gcd(p-1,4).
\]

\bigskip\noindent\textbf{Group 7.4.2}%

\[
\langle a,b,c,d\,|\,[d,a],\,[d,b],\,[d,c]\rangle
\]

The number of conjugacy classes is $p^{5}+p^{4}-p^{2}$, and the automorphism
group has order $(p-1)(p^{3}-1)(p^{3}-p)(p^{3}-p^{2})p^{15}$. The nucleus $N$
has dimension 8, and is spanned by $baa$, $bab$, $bac$, $caa$, $cab$, $cac$,
$cbb$, $cbc$. The number of orbits of subspaces of $N$ of dimension 6 under
the automorphism group is%
\[
p^{4}+p^{3}+5p^{2}+13p+(2p+3)\gcd(p-1,3)+\gcd(p-1,4)+26.
\]

The multiplier $M$ has dimension 11, and the number of orbits of the
automorphism group on subspaces of $M$ of dimension 9 is%
\begin{align*}
&  p^{6}+2p^{5}+6p^{4}+11p^{3}+27p^{2}+77p+162+(p^{2}+\frac{11}{2}p+\frac
{27}{2})\gcd(p-1,3)\\
&  +(p^{2}+\frac{1}{2}p+\frac{3}{2})\gcd(p+1,3)+(p+6)\gcd(p-1,4)+\gcd(p-1,5).
\end{align*}
However not all these subspaces are allowable. (The allowable subspaces
intersect $N$ in a subspace of dimension 6.) If a subspace of $M$ of dimension
9 is not allowable, then it either contains $N$, or intersects $N$ in a
subspace of dimension 7.

There is exactly one orbit of subspaces of $M$ of dimension 9 which contains
$N$.

There are $p+4+\gcd(p-1,3)$ orbits of subspaces of $N$ of dimension 7.
Representatives for these orbits are as follows:%

\begin{gather*}
\langle bab,bac,caa,cab,cac,cbb,cbc\rangle,\\
\langle baa,bac,caa,cab,cac+bab,cbb,cbc\rangle,\\
\langle baa,bac,caa,cab,cac+\omega bab,cbb,cbc\rangle,\\
\langle baa,bac,caa-bab,cab,cac,cbb,cbc\rangle,\\
\langle bab,bac,caa,cab,cac-baa,cbb-baa,cbc-\lambda baa\rangle,\\
\langle bab,bac,caa,cab,cac-baa,cbb-\omega baa,cbc-\lambda baa\rangle
\;(p=1\operatorname{mod}3),\\
\langle bab,bac,caa,cab,cac-baa,cbb-\omega^{2}baa,cbc-\lambda baa\rangle
\;(p=1\operatorname{mod}3),\\
\langle bab-baa,bac,caa-baa,cab,cac+baa,cbb+baa,cbc-baa\rangle,
\end{gather*}
where in the subspaces 5,6 and 7 we have $\lambda=0$, or $\lambda$ running
over a set of representatives for equivalence classes of non-zero elements in
GF$(p)$ under the equivalence relation $\lambda\sim\mu$ if $\lambda^{3}%
=\mu^{3}$.

There are 7 orbits of subspaces of $M$ of dimension 9 which intersect $N$ in
the first of these subspaces. There are 7 orbits of subspaces of $M$ of
dimension 9 which intersect $N$ in the second of these subspaces, and 4 orbits
which intersect $N$ in the third, and there are 5 orbits of subspaces of $M$
of dimension 9 which intersect $N$ in the fourth subspace.

When counting orbits of subspaces of $M$ of dimension 9 which intersect in the
fifth, sixth and seventh of these subspaces we need to distinguish between the
cases $p=1\operatorname{mod}3$ and $p=2\operatorname{mod}3$. When
$p=2\operatorname{mod}3$ we have $p$ different subspaces parametrized by
$\lambda$ with $\lambda\in\,$GF$(p)$. The number of orbits then depends on the
number of roots in GF$(p)$ of the polynomial $x^{3}-\lambda x+1$. If there are
no roots then the number is 1, if there is one root the number is 3, if there
are two roots the number is 5, and of there are three roots the number is 7.
For most integers $\lambda$ the number of roots in GF$(p)$ of $x^{3}-\lambda
x+1$ is \emph{not} PORC! But as $\lambda$ ranges over GF$(p)$ the numbers of
occurrences of 0, 1, 2 or 3 roots is predictable --- they occur $\frac{p+1}%
{3}$, $\frac{p-1}{2}$, 1, $\frac{p-5}{6}$ times respectively. When
$p=1\operatorname{mod}3$ and $\lambda=0$ there are 3 orbits of subspaces of
$M$ of dimension 9 which intersect $N$ in the fifth subspace, 1 orbit which
intersect $N$ in the sixth subspace and 1 orbit which intersect $N$ in the
seventh subspace. When $p=1\operatorname{mod}3$ and $\lambda\neq0$ we have
$p-1$ distinct subspaces of the form 5, 6 or 7. Again the number of orbits of
subspaces of $M$ of dimension 9 which intersect $N$ in one of these subspaces
depends on the number of roots of the polynomials $x^{3}-\lambda x+1$,
$x^{3}-\lambda x+\omega$, $x^{3}-\lambda x+\omega^{2}$. As above the numbers
of orbits are 1, 3, 5 or 7, with the total number of occurrences of these
numbers being $\frac{p-1}{3}$, $\frac{p-1}{2}$, 1, $\frac{p-7}{6}$.

It is perhaps worth exploring exactly what is happening here. Corresponding to
subspaces 5, 6 and 7 we have a family of algebras of dimension 7:%
\begin{gather*}
\langle a,b,c\,|\,bab,bac,caa,cab,cac-baa,cbb-baa,cbc-\lambda baa,\text{ class
3}\rangle,\\
\langle a,b,c\,|\,bab,bac,caa,cab,cac-baa,cbb-\omega baa,cbc-\lambda
baa,\text{ class 3}\rangle,\\
\langle a,b,c\,|\,bab,bac,caa,cab,cac-baa,cbb-\omega^{2}baa,cbc-\lambda
baa,\text{ class 3}\rangle.
\end{gather*}
(The last two of these presentations only arise when $p=1\operatorname{mod}%
3$.)The automorphism groups of these algebras depend on the numbers of roots
of the polynomials $x^{3}-\lambda x+1$, $x^{3}-\lambda x+\omega$,
$x^{3}-\lambda x+\omega^{2}$. When $p=2\operatorname{mod}3$ we only need to
consider the first of these three presentations, and the automorphism group
has order $(p^{3}-1)p^{12}$, $(p-1)(p^{2}-1)p^{12}$, $(p-1)^{2}p^{13}$,
$(p-1)^{3}p^{12} $ depending on whether $x^{3}-\lambda x+1$ has 0, 1, 2 or 3
roots. When $p=1\operatorname{mod}3$ and $\lambda\neq0$ then the automorphism
groups have the same orders as in the case $p=2\operatorname{mod}3$, but now
depending on the numbers of roots of $x^{3}-\lambda x+1$ for the first
presentation, the number of roots of $x^{3}-\lambda x+\omega$ for the second
presentation, and depending on the number of roots of $x^{3}-\lambda
x+\omega^{2}$ for the third presentation. When $x=0$ and
$p=1\operatorname{mod}3$ then the automorphism groups have orders
$3(p-1)^{2}(p^{2}-1)p^{12}$, $3(p^{3}-1)p^{12}$, $3(p^{3}-1)p^{12}$.

There are three orbits of subspaces of $M$ of dimension 9 which intersect $N$
in the last of these subspaces.

So the total number of orbits of non-allowable subalgebras of $M$ of dimension 
9 is $3p+24+\gcd(p-1,3)$. The number of immediate descendants of
dimension 9 can be obtained by subtracting away this number from the total 
number of orbits of subspaces of $M$ dimension 9, and so is
\begin{align*}
&  p^{6}+2p^{5}+6p^{4}+11p^{3}+27p^{2}+74p+138+(p^{2}+\frac{11}{2}p+\frac
{25}{2})\gcd(p-1,3)\\
&  +(p^{2}+\frac{1}{2}p+\frac{3}{2})\gcd(p+1,3)+(p+6)\gcd(p-1,4)+\gcd(p-1,5).
\end{align*}

\textbf{Group 7.4.3}%
\[
\langle a,b,c,d\,|\,[d,a],\,[c,b],\,[d,c]\rangle
\]

The number of conjugacy classes is $3p^{4}-3p^{2}+p$, and the automorphism
group has order $2(p-1)^{4}p^{16}$. The nucleus has dimension 8, and is
spanned by $baa$, $bab$, $bac$, $bad$, $caa$, $cac$, $dbb$, $dbd$. The number
of orbits of subspaces of dimension 6 in the nucleus under the automorphism
group is%
\begin{align*}
&  \frac{1}{2}(p^{5}+4p^{4}+15p^{3}+48p^{2}+139p+477+\gcd(p-1,3)(4p+26)\\
&  +\gcd(p-1,4)(3p+15)+2\gcd(p-1,5)).
\end{align*}

The multiplier $M$ has dimension 11, and the number of orbits of subspaces of
$M$ of dimension 9 under the automorphism group is%
\begin{align*}
&  \frac{1}{2}(p^{5}+4p^{4}+21p^{3}+85p^{2}+355p+1875)\\
&  +\frac{1}{2}(p^{2}+13p+103)\gcd(p-1,3)\\
&  +\frac{1}{2}(3p+44)\gcd(p-1,4)\\
&  +\gcd(p-1,3)\gcd(p-1,4)+3\gcd(p-1,5).
\end{align*}
The number of descendants of dimension 9 is the number of orbits of allowable
subspaces of $M$ of dimension 9. So to we need to find the number of orbits of
non-allowable subspaces of $M$ of dimension 9. A non-allowable subspace of $M$
of dimension 9 must either contain the nucleus, or intersects the nucleus in a
subspace of dimension 7. To compute these subspaces we adopt a different
presentation for $L$.%

\[
\langle a,b,c,d\,|\,ca,\,da,\,db,\,\text{class 2}\rangle.
\]

With this presentation the nucleus of dimension 8 is spanned by $baa$, $bab$,
$bac$, $cbb$, $cbc$, $cbd$, $dcc$, $dcd$.

The automorphism group is determined by its action on the Frattini quotient of
the group and is given by matrices in GL$(4,p)$ of the following two forms:%
\[
\left(
\begin{array}
[c]{cccc}%
\alpha & 0 & 0 & 0\\
\beta & \gamma & 0 & \delta\\
\varepsilon & 0 & \zeta & \eta\\
0 & 0 & 0 & \theta
\end{array}
\right)  ,
\]%
\[
\left(
\begin{array}
[c]{cccc}%
0 & 0 & 0 & \alpha\\
\beta & 0 & \gamma & \delta\\
\varepsilon & \zeta & 0 & \eta\\
\theta & 0 & 0 & 0
\end{array}
\right)  .
\]

We consider the action of these automorphisms on the dual space of the
nucleus. It turns out that there are $p+26$ orbits of one dimensional
subspaces of this dual space. A set of basis elements for representatives of
these orbits is as follows:%
\begin{align*}
&  (0,0,0,0,0,0,0,1),\\
&  (0,0,0,0,0,0,1,0),\\
&  (0,0,0,0,0,1,0,0),\\
&  (0,0,0,0,0,1,0,1),\\
&  (0,0,0,0,0,1,1,0),\\
&  (0,0,0,0,0,1,1,1),\\
&  (0,0,0,0,1,0,0,0),\\
&  (0,0,0,0,1,0,0,1),\\
&  (0,0,0,1,0,0,0,k),
\end{align*}
where $k$ is not a square,%
\begin{align*}
&  (0,0,0,1,0,0,1,0),\\
&  (0,0,0,1,0,1,1,k),
\end{align*}
where $1+k$ is not a square,%
\begin{align*}
&  (0,0,0,1,1,0,0,0),\\
&  (0,0,1,0,0,0,0,1),\\
&  (0,0,1,0,0,0,1,0),\\
&  (0,0,1,0,0,1,0,0),\\
&  (0,0,1,0,0,1,1,0),\\
&  (0,1,0,0,0,0,0,1),\\
&  (0,1,0,0,0,0,1,0),\\
&  (0,1,0,0,0,1,1,0),\\
&  (0,1,0,0,0,1,1,1),\\
&  (0,1,1,0,0,0,0,1),\\
&  (0,1,1,0,0,1,\lambda,0)\,,
\end{align*}

with $1\leq\lambda\leq p-1$, and%
\begin{align*}
&  (1,0,0,0,0,0,0,1),\\
&  (1,0,0,0,0,1,0,1),\\
&  (1,0,0,0,0,1,1,-1),\\
&  (1,0,0,0,0,1,1,-k),\\
&  (1,0,0,0,k,0,0,1),
\end{align*}
where in the last two $k$ is not a square, and finally%
\[
(1,0,0,1,0,1,1,k)
\]
where neither $1+k$ nor $-k$ are squares. The annihilators of these elements
of the dual space of $N$ give subspaces of $N$ of dimension 7.

There are 3 orbits of subspaces of $M$ of dimension 9 which contain the
nucleus. For the first 21 of the subspaces $S$ giving representatives for the
$p+26$ orbits of subspaces of the nucleus $N$ of dimension 7, the number of
orbits of subspaces of $M$ of dimension 9 which contain subspaces intersecting
$N$ in $S$ is as follows: 7, 9, 7, 7, 5, $\frac{p+7}{2}$, 9, 5, 7, 7,
$\frac{p+7}{2}$, 5, 9, 8, 3, $p+2$, 7, $4+\gcd(p-1,3)$, $p+2$, $4+\gcd
(p-1,5)$, $p+6$. For the subspaces corresponding to $(0,1,1,0,0,1,\lambda,0)$,
the relevant number of orbits is $\frac{p+3}{2}$ if $-\lambda$ is not a square
in GF$(p)$, $\frac{p+5}{2}$ if $-\lambda$ is a square in GF$(p)$ and
$\lambda\neq-4$, and $\frac{p+7}{2}$ if $\lambda=-4$. For the last six
subspaces the relevant numbers of orbits are
\[
4+\gcd(p-1,3),
\]%
\[
\frac{1}{2}(p+7+2\gcd(p-1,3)+\gcd(p-1,4)),
\]%
\[
\frac{1}{8}(p^2+6p+9+2\gcd(p-1,4)),
\]%
\[
\frac{p^{2}-1}{4}+p+2,
\]%
\[
\frac{p-1}{2}+4+\frac{1}{2}\gcd(p-1,4)+\gcd(p-1,3),
\]%
\[
\frac{1}{8}(p^2+6p+9+2\gcd(p-1,4)).
\]

So the total number of non-allowable subspaces of $M$ of dimension 9 is
\[
p^2+9p+137+4\gcd(p-1,3)+\frac{3}{2}\gcd(p-1,4)+\gcd(p-1,5).
\]
Subtracting this total from the number of orbits of subspaces of $M$ of 
dimension 9 we obtain the number of descendants of dimension 9:
\begin{align*}
&  \frac{1}{2}(p^{5}+4p^{4}+21p^{3}+83p^{2}+337p+1601)\\
&  +\frac{1}{2}(p^{2}+13p+95)\gcd(p-1,3)\\
&  +\frac{1}{2}(3p+41)\gcd(p-1,4)\\
&  +\gcd(p-1,3)\gcd(p-1,4)+2\gcd(p-1,5).
\end{align*}

\bigskip\noindent\textbf{Group 7.4.4}%

\[
\langle a,b,c,d\,|\,[d,a]=[c,b],\,[d,b],\,[d,c]\rangle
\]

The number of conjugacy classes is $2p^{4}+p^{3}-2p^{2}$, and the automorphism
group has order $(p-1)(p^{2}-1)(p^{2}-p)p^{17}$. The nucleus has dimension 8,
and is spanned by $baa$, $bab$, $bac$, $bad$, $caa$, $cab$, $cac$, $cad$. The
number of orbits of subspaces of dimension 6 in the nucleus under the
automorphism group is%
\[
p^{3}+3p^{2}+25p+64+\gcd(p-1,3)+2\gcd(p-1,4).
\]

The multiplier has dimension 11, and the number of orbits of subspaces of the
multiplier of dimension 9 is%
\[
p^{3}+5p^{2}+52p+222+2\gcd(p-1,3)+3\gcd(p-1,4)
\]
when $p>5$, and 748 when $p=5$. Not all these subspaces are allowable. A
subspace which is not allowable must either contain the nucleus, or must
intersect the nucleus in a subspace of dimension 7. There are 2 orbits of 9
dimensional subspaces of the multiplier which contain the nucleus. There are
$p+10$ orbits of subspaces of dimension 7 in the nucleus. Representatives for
these subspaces are given by the following one dimensional subspaces of the
dual space of the nucleus:%
\begin{gather*}
\langle(0,0,0,0,0,0,0,1)\rangle,\\
\langle(0,0,0,0,0,0,1,0)\rangle,\\
\langle(0,0,0,0,0,1,0,0)\rangle,\\
\langle(0,0,0,0,1,0,0,0)\rangle,\\
\langle(0,0,0,1,0,0,1,0)\rangle,\\
\langle(0,0,0,1,0,1,0,0)\rangle,\\
\langle(0,0,0,1,-\omega,0,0,0)\rangle,\\
\langle(0,0,1,0,0,x,0,0)\rangle\;x\neq0,\;x\symbol{126}x^{-1},\\
\langle(0,0,1,0,0,-1,1,0)\rangle,\\
\langle(0,1,0,0,0,0,-\omega,0)\rangle,\\
\langle(0,1,0,0,0,1,x,0)\rangle\;\text{all }x\text{ with }1-4x\text{ not
square},\\
\langle(0,1,0,0,1,0,0,0)\rangle.
\end{gather*}
(Here $x\symbol{126}x^{-1}$ indicates that $x$ and $x^{-1}$ give subspaces in
the same orbit.) For the first seven of these subspaces $S$, there are
(respectively) 4, 5, 6, 7, 5, $\frac{p+5}{2}$, $\frac{p+5}{2}$ orbits of
subspaces of $M$ of dimension 9 which contain representatives intersecting $N
$ in $S$. The next subspace $S$ is indexed by a parameter $x$, which runs
through a set of representatives for equivalence classes of non-zero elements
of GF$(p)$ under the equivalence relation $x\symbol{126}x^{-1}$. When $x=1$ we
have 3 orbits of subspaces of $M$ of dimension 9 which contain representatives
intersecting $N$ in $S$. When $x=-1$ we have 3 orbits when $p=5$ and 2 orbits
when $p>5$. When $x\neq\pm1$ we have 4 orbits, except (for $p>5$) when
$x=\frac{2}{3}$ or $\frac{3}{2}$ we have 5 orbits. The next subspace gives 4
orbits if $p=5$ and 3 orbits if $p>5$. The next subspace gives 2 orbits. The
penultimate subspace is indexed by $x$, where $x$ takes on the $\frac{p-1}{2}$
values in GF$(p)$ such that $1-4x$ is not square. Each of the corresponding
subspaces gives 2 orbits. The last subspace gives 5 orbits.

So when $p=5$ there are 64 orbits of non-allowable subspaces of $M$ of 
dimension 9, and when $p>5$ there are $4p+43$ orbits. By subtracting
these numbers from the grand total of all orbits of 9 dimensional subspaces of
$M$, we obtain the number of descendants of dimension 9. When $p=5$ there are
684 descendants, and when $p>5$ the number of descendants is
\[
p^{3}+5p^{2}+48p+179+2\gcd(p-1,3)+3\gcd(p-1,4).
\]

\bigskip\noindent\textbf{Group 7.4.5}%

\[
\langle a,b,c,d\,|\,[c,a],\,[d,a]=[c,b],\,[d,b]\rangle
\]

The number of conjugacy classes is $2p^{4}+p^{3}-2p^{2}$, and the automorphism
group has order $(p^{2}-1)^{2}(p^{2}-p)p^{13}$. The nucleus has dimension 8
and is spanned by $baa$, $bab$, $bac$, $bad$, $cbc$, $cbd$, $dcc$, $dcd$. The
number of orbits of subspaces of dimension 6 in the nucleus under the
automorphism group is%
\begin{align*}
&  p^{6}+p^{5}+5p^{4}+6p^{3}+18p^{2}+22p+51+(4p+11)\gcd(p-1,3)\\
&  +2\gcd(p-1,4)+\gcd(p-1,5).
\end{align*}

The multiplier $M$ has dimension 11, and the number of orbits under the
automorphism group of subspaces of $M$ of dimension 9 is%
\begin{align*}
&  p^{6}+p^{5}+6p^{4}+13p^{3}+34p^{2}+72p+205+(p^{2}+6p+39)\gcd(p-1,3)\\
&  +4\gcd(p-1,4)+\gcd(p-1,5)+\gcd(p-1,9).
\end{align*}
The number of descendants of dimension 9 is the number of orbits of allowable
subspaces of $M$ of dimension 9. So to compute the number of descendants of
dimension 9 we need to compute the number of orbits of non-allowable subspaces
of $M$ of dimension 9. A non-allowable subspace of $M$ of dimension 9 must
either contain the nucleus, or intersect the nucleus in a subspace of
dimension 7. To compute these subspaces we adopts a different presentation for
$L$.%

\[
\langle a,b,c,d\,|\,ba,\,da-cb,\,dc,\,\text{class 2}\rangle.
\]
With this presentation, nucleus of dimension 8, spanned by $caa$, $cab$, $cac
$, $cad$, $cbb$, $cbd$, $dbb$, $dbd$. The action of the automorphism group on
$L/L^{2}$ is given by matrices in GL$(4,p)$ of the form%
\[
\left(
\begin{array}
[c]{cccc}%
\lambda\alpha & \lambda\beta & \mu\alpha & \mu\beta\\
\lambda\gamma & \lambda\delta & \mu\gamma & \mu\delta\\
\nu\alpha & \nu\beta & \xi\alpha & \xi\beta\\
\nu\gamma & \nu\delta & \xi\gamma & \xi\delta
\end{array}
\right)  ,
\]
with $(\alpha\delta-\beta\gamma)(\lambda\xi-\mu\nu)\neq0$

As usual we compute the orbits of one dimensional subspaces of the dual space
of the nucleus under the action of the automorphism group. It turns out that
there are $2p+10+\gcd(p-1,3)$ orbits of these subspaces. Bases for
representative subspaces are as follows:%
\begin{align*}
&  (0,0,0,0,0,0,0,1),\\
&  (0,0,0,0,0,1,0,0),\\
&  (0,0,0,0,0,1,1,0),\\
&  (0,0,0,1,0,0,0,1),\\
&  (0,0,0,1,0,0,0,\omega),\\
&  (0,0,0,1,0,0,1,0),\\
&  (0,0,0,1,1,0,0,0),\\
&  (0,0,0,1,1,0,0,1),\\
&  (0,0,0,1,1,0,0,\omega),
\end{align*}

\[
(0,0,1,0,0,0,0,\omega)\text{ if }p=1\operatorname{mod}3,
\]

\[
(0,0,1,0,0,0,1,0),
\]

\[
(0,0,1,0,0,1,0,\lambda)\text{ if }p=2\operatorname{mod}3\text{ for any one
}\lambda\text{ such that }x^{3}+3x+\lambda\text{ is irreducible,}%
\]

\begin{align*}
&  (0,0,1,0,0,1,1,0),\\
&  (0,0,1,0,0,\omega,1,0),\\
&  (0,0,1,0,1,0,0,2),
\end{align*}

\begin{align*}
(0,0,1,0,1,0,0,2\omega)\text{ if }p  &  =1\operatorname{mod}3,\\
(0,0,1,0,1,0,0,2\omega^{2})\text{ if }p  &  =1\operatorname{mod}3,
\end{align*}

\[
(0,0,1,0,1,0,1,\lambda)\text{ for all }\lambda\neq0,\frac{1}{27},
\]

\[
(0,1,1,0,0,\lambda,\mu,\nu)\text{, (}p-1\text{ of them).}%
\]
For each of these elements $f$ of the dual space of the nucleus we take a set
of basis elements of the annihilator of $f$ as relators.

There are 3 orbits of subspaces of $M$ of dimension 9 which contain the
nucleus. For the first 14 of the $2p+10+\gcd(p-1,3)$ subspaces $S$ giving
representatives for the $2p+10+\gcd(p-1,3)$ orbits of subspaces of the nucleus
$N$ of dimension 7, the number of orbits of subspaces of $M$ of dimension 9
which contain subspaces intersecting $N$ in $S$ is as follows: 7, 7, 5, 7, 7,
$2p+5+\gcd(p-1,4)$, $p+4$, $\frac{1}{2}(p^{2}+2p+3)$, $\frac{1}{2}%
(p^{2}+2p+3)$, 7, $\frac{1}{2}(p+3)(\gcd(p-1,3)+1)+1$, 7, $\frac{1}{2}%
(p^{2}+2p+3)$, $\frac{1}{2}(p^{2}+2p+3)$. For the next subspace the relevant
number of orbits is $\frac{1}{2}(p^{2}+2p+3)$ when $p=2\operatorname{mod}3$,
and $\frac{1}{12}(p^{2}+4p+31)$ when $p=1\operatorname{mod}3$. The next two
subspaces only arise when $p=1\operatorname{mod}3$, and then the relevant
number of orbits is $\frac{1}{3}(p^{2}+p+7)\,$\ in both cases. The next
subspace is actually a family of $p-2$ subspaces indexed by a parameter
$\lambda$. It turns out that the relevant number of orbits depends on the
number of roots over GF$(p) $ of the polynomial $2x^{3}+x^{2}-\lambda$. The
discriminant of this polynomial is $\frac{1}{4}\lambda(1-27\lambda)$, and the
discriminant is zero if $\lambda=0$ or $\frac{1}{27}$, the two excluded values
of $\lambda$. So the number of roots for $\lambda\neq0,\frac{1}{27}$ is 0, 1
or 3. If there are $n$ roots, then the relevant number of orbits is then%
\[
\frac{p^{2}+p+1+n(p+2)}{n+1}.
\]
Intriguingly, if we pick an integer $\lambda$ such that the Galois group of
$2x^{3}+x^{2}-\lambda$ over the rationals is $S_{3}$, then the number of roots
of $2x^{3}+x^{2}-\lambda$ over GF$(p)$ is not PORC. But as $\lambda$ ranges
over GF$(p)$ the number of times 0, 1, 2, or 3 roots arises is PORC. The
following table gives the number of times $n$ roots arises $(n=0,1,2,3)$ as
$\lambda$ ranges over GF$(p)$.%
\[%
\begin{tabular}
[c]{|c|c|c|c|c|}\hline
$n$ & 0 & 1 & 2 & 3\\\hline
$p=1\operatorname{mod}3$ & $\frac{p-1}{3}$ & $\frac{p-1}{2}$ & 2 & $\frac
{p-7}{6}$\\\hline
$p=2\operatorname{mod}3$ & $\frac{p+1}{3}$ & $\frac{p-3}{2}$ & 2 & $\frac
{p-5}{6}$\\\hline
\end{tabular}
\]
We obtain the figures for $n=0$ by calculating the number of irreducible
polynomials of the form $2x^{3}+x^{2}-\lambda$, and the figures for $n=1$ by
calculating the number of values of $\lambda$ for which the discriminant is
not a square. We get 2 roots if $\lambda=0$ or $\frac{1}{27}$, and the figure
for $n=3$ then follows from the fact that the rows have to add up to $p$. So
the sum of all the relevant number of orbits over these subspaces is PORC.

Finally we have $p-1$ subspaces corresponding to elements in the dual space of
the nucleus of the form $(0,1,1,0,0,\lambda,\mu,\nu)$ with $\lambda,\mu\neq0$.
For $\frac{p-1}{2}$ of these subspaces the relevant number of orbits is
$\frac{1}{2}(p^{2}+2p+3)$, and for $\frac{p-5}{2}$ of the subspaces the
relevant number of orbits is $\frac{1}{4}(p^{2}+4p+7)$. Finally we have the
subspaces corresponding to%
\[
(0,1,1,0,0,\omega,\omega,0)\text{ and }(0,1,1,0,0,-3\omega,-27\omega,0).
\]
For the first of these the relevant number of orbits is $p+2$ when
$p=1\operatorname{mod}3$ and $2p+3$ when $p=2\operatorname{mod}3$, and for the
second the relevant number of orbits is $p+2$.

So the total number of non-allowable subspaces of $M$ of dimension 9 is
\[
p^3+3p^2+10p+56+4\gcd(p-1,3)+\gcd(p-1,4)
\]
We obtain the number of descendants of dimension 9 by subtracting this
number from the total number of subspaces of $M$ of dimension 9:
\begin{align*}
&  p^{6}+p^{5}+6p^{4}+12p^{3}+31p^{2}+62p+149+(p^{2}+6p+35)\gcd(p-1,3)\\
&  +3\gcd(p-1,4)+\gcd(p-1,5)+\gcd(p-1,9).
\end{align*}

\bigskip\noindent\textbf{Group 7.4.6}%

\[
\langle a,b,c,d\,|\,[d,b]=[c,a]^{\omega},\,[d,c]=[b,a],\,[c,b]\rangle
\]

The number of conjugacy classes is $p^{4}+2p^{3}-p^{2}-p$, and the
automorphism group has order $2(p^{2}-1)^{2}p^{16}$. The nucleus has dimension
8 and is spanned by $baa$, $bac$, $bad$, $caa$, $cac$, $cad$, $daa$, $dad$.
The number of orbits of subspaces of dimension 6 in the nucleus under the
automorphism group is%
\[
\frac{1}{2}(p^{5}+5p^{3}+8p^{2}+23p+29-2p\gcd(p-1,3)+(p+1)\gcd(p-1,4)).
\]
The $p$-multiplier $M$ has dimension 11, and the number of orbits of subspaces
of codimension 2 in $M$ under the automorphism group is%
\[
\frac{1}{2}(p^{5}+5p^{3}+17p^{2}+45p+137-(p^{2}+3p+13)\gcd(p-1,3)+(p+4)\gcd
(p-1,4)).
\]

Not all subspaces of $M$ of dimension 9 are allowable, and to obtain the total
number of descendants of 7.4.6 of order $p^{9}$ with exponent $p$ we need to
subtract the number of orbits of non-allowable subspaces of dimension 9 from
the total number of orbits of subspaces of dimension 9. A non-allowable
subspace of dimension 9 either contains the nucleus, or intersects the nucleus
in a subspace of dimension 7. There are 2 orbits of subspaces of $M$ of
dimension 9 containing the nucleus.

To compute the number of orbits of subspaces of $M$ of dimension 9 which
intersect the nucleus in a subspace of dimension 7, we switch to an
alternative presentation for group 7.4.6:%
\[
\langle a,b,c,d\,|\,da,\,db-\omega ca,\,dc-ba,\,\text{class 2}\rangle.
\]
With this presentation the nucleus of dimension 8, spanned by $baa$, $bab$,
$bac$, $caa$, $cab$, $cac$, $cbb$, $cbc$.

The action of the automorphism group on $L/L^{2}$ is given by the following
matrices in GL$(4,p)$:%
\[
\left(
\begin{array}
[c]{cccc}%
\alpha & 0 & 0 & \beta\\
\lambda & \delta & \omega\gamma & \mu\\
\nu & \pm\gamma & \pm\delta & \xi\\
\pm\omega\beta & 0 & 0 & \pm\alpha
\end{array}
\right)  .
\]

As usual we compute the orbits of one dimensional subspaces of the dual space
of the nucleus under the action of the automorphism group. It turns out that
there are $p+5$ orbits of these subspaces. Bases for representative subspaces
are as follows:%
\begin{align*}
&  (0,0,0,0,0,0,0,1),\\
&  (0,0,0,0,0,1,0,0),\\
&  (0,0,0,1,0,0,0,0),\\
&  (0,0,0,1,0,0,0,1),
\end{align*}%
\[
(0,0,0,1,0,0,1,\lambda)\text{ for one }\lambda\text{ such that }\omega
^{2}\lambda^{2}-\omega\text{ is not a square,}%
\]
and $p$ elements%
\[
(0,0,1,0,\lambda,\mu,0,0)
\]
with $\lambda\neq0$. The numbers of orbits of subspaces of $M$ of dimension 9
corresponding to the first five of these are 3,
\begin{align*}
&  3,\\
&  \frac{1}{2}(p+5),\\
&  6-\gcd(p-1,3),\\
&  \frac{1}{4}(p^{2}+(6-\gcd(p-1,4))p+5+\gcd(p-1,4)),\\
&  \frac{1}{4}(p^{2}+\gcd(p-1,4)p+3+\gcd(p-1,4)).
\end{align*}
In the last case, the number of orbits of subspaces of $M$ of dimension 9
corresponding to the $p$ representative subspaces are as follows. We have
three representative subspaces giving $6-\gcd(p-1,3)$, $2$ and $\frac{p+7}{2}$
orbits respectively. We have $\frac{p-1}{2}$ representative subspaces each
giving $\frac{p+3}{2}$ orbits. And finally we have $\frac{p-5}{2}$
representative subspaces each giving $\frac{p+5}{2}$ orbits.

So the total number of orbits of non-allowable subspaces of $M$ of dimension 9 is
\[
p^2+3p+20-2\gcd(p-1,3)+\frac{1}{2}\gcd(p-1,4).
\]
We obtain the number of descendants of dimension 9 by subtacting this
number from the total number of orbits of subspaces of $M$ of dimension
9, giving:
\[
\frac{1}{2}(p^{5}+5p^{3}+15p^{2}+39p+97-(p^{2}+3p+9)\gcd(p-1,3)+(p+3)\gcd(p-1,4)).
\]

\subsection{Five generator groups}

\noindent\textbf{Group 7.5.1}%

\[
\langle
a,b,c,d,e\,|\,[c,b],\,[d,a],\,[d,b],\,[d,c],\,[e,a],\,[e,b],\,[e,c],\,[e,d]\rangle
\]

The number of conjugacy classes is $2p^{5}-p^{3}$ and the order of the
automorphism group is $(p-1)(p^{2}-1)^{2}(p^{2}-p)^{2}p^{18}$. The nucleus has
dimension 5, and is generated by $baa$, $bab$, $bac$, $caa$, $cac$. From the
calculation of the descendants of $\langle a,b,c\,|\,[c,b],$ class 2, exponent
$p\rangle$, we may assume that one of the following holds:%

\[
caa=cab=cac=0,
\]%
\[
bab=bac=cab=cac=0,
\]%
\[
bac=cac=0,\,caa=bab,
\]%
\[
baa=bac=cac=0,
\]%
\[
bac=caa=0,\,cac=bab,
\]%
\[
bac=caa=0,\,cac=\omega bab,
\]

\[
bac=0,\,caa=baa,\,cac=-bab,
\]%
\[
baa=bac=caa=0,
\]%
\[
bac=caa=0,\,baa=cac,
\]%
\[
bac=0,\,baa=cac,\,caa=bab,
\]%
\[
bac=0,\,baa=cac,\,caa=\omega bab,\,(p=1\operatorname{mod}3)
\]%
\[
baa=caa=cac=0,
\]%
\[
baa=cac=0,\,caa=bab,
\]%
\[
caa=cac=0,\,baa=bab,
\]%
\[
baa=caa=0,\,cac=\omega bab,
\]%
\[
baa=0,\,caa=bac,\,cac=\omega bab,
\]%
\[
baa=0,\,caa=kbab+bac,\,cac=\omega bab,\,(p=2\operatorname{mod}3),
\]
where $k$ is any element of $\mathbb{Z}_{p}$ which is not a value of
\[
\frac{\lambda(\lambda^{2}+3\omega\mu^{2})}{\mu(3\lambda^{2}+\omega\mu^{2})}.
\]

The number of descendants of class 3 with order $p^{9}$ and exponent $p$ in
each of the first 9 cases above is%
\begin{gather*}
2p+60+2\gcd(p-1,3),\\
4p+27,\\
3p+25,\\
p+61+\gcd(p-1,4),\\
\frac{1}{2}(p^{2}+11p+56+3\gcd(p-1,4))/2,\\
\frac{1}{2}(p^{2}+11p+56+3\gcd(p-1,4))/2,\\
p+49,\\
\frac{1}{2}(p^{2}+6p+67+9\gcd(p-1,3)+\gcd(p-1,5),\\
p^{3}+3p^{2}+9p+27+3\gcd(p-1,3)+3\gcd(p-1,4).
\end{gather*}
When $p=1\operatorname{mod}3$ then the number of descendants of class 3 with
order $p^{9}$ and exponent $p$ in cases 10 and 11 above combined is%
\[
\frac{1}{2}(p^{4}+2p^{3}+4p^{2}+14p+39)
\]
and when $p=2\operatorname{mod}3$ the number is%
\[
\frac{1}{2}(p^{4}+2p^{3}+4p^{2}+8p+15).
\]
The number of descendants of class 3 with order $p^{9}$ and exponent $p$ in
each of cases 12, 13, 14, 15 above is%

\begin{gather*}
2p+35,\\
p^{3}+3p^{2}+7p+27+3\gcd(p-1,3),\\
p^{2}+7p+21,\\
\frac{1}{2}(p^{2}+4p+33-\gcd(p-1,3)+\gcd(p+1,5)).
\end{gather*}
Finally, in cases 16 and 17 above combined, the number of descendants is%
\[
\frac{1}{2}(p^{4}+2p^{3}+4p^{2}+8p+15)
\]
when $p=1\operatorname{mod}3$, and
\[
\frac{1}{2}(p^{4}+2p^{3}+4p^{2}+10p+23)
\]
when $p=2\operatorname{mod}3$.

\bigskip\noindent\textbf{Group 7.5.2}%

\[
\langle
a,b,c,d,e\,|\,[c,b],\,[d,a],\,[d,b]=[b,a],\,[d,c],\,[e,a],\,[e,b],\,[e,c],\,[e,d]\rangle
\]

The number of conjugacy classes is $p^{5}+2p^{4}-p^{3}-2p^{2}+p$ and the order
of the automorphism group is $2(p-1)(p^{2}-1)^{2}(p^{2}-p)^{2}p^{14}$. The
nucleus has dimension 4. It is convenient to take another presentation for
this group:%
\[
\langle
a,b,c,d,e\,|\,[c,a],\,[c,b],\,[d,a],\,[d,b],\,[e,a],\,[e,b],\,[e,c],\,[e,d]\rangle
.
\]

The nucleus has dimension 4 and is spanned by $baa$, $bab$, $dcc$, $dcd$. We
need to factor out a two dimensional subspace of the nucleus, and we can
assume that one of the following pairs of relations holds:%
\begin{align*}
dcc  &  =dcd=0,\\
bab  &  =dcd=0,\\
bab-dcc  &  =dcd=0,\\
dcc-baa  &  =dcd-bab=0,
\end{align*}
or, when $p=1\operatorname{mod}3$,%
\[
dcc-baa=dcd-\omega bab=0\text{.}%
\]

The number of descendants of order $p^{9}$ with class 3 and exponent $p$ in
the first 3 cases is%
\begin{align*}
&  p+65,\\
&  \frac{1}{2}(4p+135+5\gcd(p-1,3)),\\
&  6p+93+5\gcd(p-1,3)+4\gcd(p-1,4),
\end{align*}
and the number of descendants in the last two cases combined is%
\[
\frac{1}{2}(2p^{2}+13p+11+30\gcd(p-1,3)+\gcd(p-1,4)).
\]

\bigskip\noindent\textbf{Group 7.5.3}%

\[
\langle
a,b,c,d,e\,|\,[c,b],\,[d,a],\,[d,b]=[c,a],\,[d,c],\,[e,a],\,[e,b],\,[e,c],\,[e,d]\rangle
\]

The number of conjugacy classes is $p^{5}+p^{4}-p^{2}$ and the order of the
automorphism group is $(p-1)^{2}(p^{2}-1)(p^{2}-p)p^{18}$. The nucleus has
dimension 4, and is spanned by $baa$, $bab$, $bac$, $bad$. We may assume that
one of the following four sets of relations holds:%
\begin{align*}
bac  &  =bad=0,\\
baa  &  =bac=0,\\
bab  &  =bac=0,\\
baa  &  =bab=0.
\end{align*}

The number of descendants of order $p^{9}$ with class 3 and exponent $p$ in
these four cases are:%
\begin{align*}
&  4p+48+\gcd(p-1,3),\\
&  3p+69,\\
&  5p+46,\\
&  2p+44.
\end{align*}

\bigskip\noindent\textbf{Group 7.5.4}%

\[
\langle a,b,c,d,e\,|\,[c,b],\,[d,a],\,[d,b]=[c,a],\,[d,c]=[b,a]^{\omega
},\,[e,a],\,[e,b],\,[e,c],\,[e,d]\rangle
\]

The number of conjugacy classes is $p^{5}+p^{3}-p$ and the order of the
automorphism group is $2(p-1)(p^{4}-1)(p^{4}-p^{2})p^{14}$. The nucleus has
dimension 4, and is spanned by $baa$, $bab$, $bac$, $bad$. We may assume that
one of the following three sets of relations holds:%
\begin{align*}
bac  &  =bad=0,\\
baa  &  =bad=0,
\end{align*}
or when $p=2\operatorname{mod}3$%
\[
bac=bad-xbaa=0,
\]
where $x$ is chosen so that $y^{3}-y-\frac{4x}{3(x^{2}+3)}$ is irreducible.
Note that $x^{2}+3$ can never be zero since $p=2\operatorname{mod}3$. Also
note that $x$ cannot be zero. We also require that $x^{2}\neq3$.

The number of descendants of order $p^{9}$ with class 3 and exponent $p$ in
the first and last cases combined is%
\[
\frac{1}{2}(2p^{2}+13p+65-8\gcd(p-1,3)+\gcd(p-1,4)).
\]
The number of descendants in the second case is%
\[
\frac{1}{2}(4p+35-3\gcd(p-1,3)).
\]

\bigskip\noindent\textbf{Group 7.5.5}%

\[
\langle
a,b,c,d,e\,|\,[c,b],\,[d,a],\,[d,b],\,[d,c],\,[e,a],\,[e,b],\,[e,c],\,[e,d]=[b,a]\rangle
\]

The number of conjugacy classes is $p^{5}+p^{4}-p^{2}$ and the order of the
automorphism group is $(p-1)^{2}(p^{2}-1)(p^{2}-p)p^{15}$. The nucleus has
dimension 2. The number of descendants of order $p^{9}$ with exponent $p$ is%
\[
p^{3}+4p^{2}+16p+4\gcd(p-1,4)+6\gcd(p-1,3)+80.
\]

\bigskip\noindent\textbf{Group 7.5.6}%

\[
\langle
a,b,c,d,e\,|\,[c,b],\,[d,a],\,[d,b]=[c,a],\,[d,c],\,[e,a],\,[e,b],\,[e,c],\,[e,d]=[b,a]\rangle
\]

The number of conjugacy classes is $p^{5}+p^{3}-p$ and the order of the
automorphism group is $(p-1)(p^{2}-1)(p^{2}-p)p^{14}$. This group has no class
3 descendants of order $p^{9}$.

\subsection{Six generator groups}

\noindent\textbf{Group 7.6.1}%

\[
\langle a,b\rangle\times\langle c\rangle\times\langle d\rangle\times\langle
e\rangle\times\langle f\rangle
\]

The number of conjugacy classes is $p^{6}+p^{5}-p^{4}$, and the automorphism
group has order $(p^{7}-p^{5})(p^{7}-p^{6})(p^{5}-p)(p^{5}-p^{2})(p^{5}%
-p^{3})(p^{5}-p^{4})$. The nucleus has dimension 2. If $G$ is an immediate
descendant of 7.6.1 of order $p^{9}$ with exponent $p$, then the subgroup of
$G$ generated by $a,b$ must be the free class 3 group with exponent $p$ of
rank 2, and has order $p^{5}$. The subgroup of $G$ generated by $c,d,e,f$ is
either elementary abelian, or is one of 5.4.1, 5.4.2, 6.4.1, 6.4.2, 6.4.3,
6.4.4. The number of descendants in the seven cases is 8, 15, 2, 15, 3, 3, 2,
making 48 descendants in all.

\bigskip\noindent\textbf{Group 7.6.2}%

\[
\langle a,b\rangle\times_{\lbrack b,a]=[d,c]}\langle c,d\rangle\times\langle
e\rangle\times\langle f\rangle
\]

The number of conjugacy classes is $p^{6}+p^{3}-p^{2}$, and the automorphism
group has order $(p^{7}-p^{3})(p^{7}-p^{6})(p^{5}-p^{3})p^{4}(p^{3}%
-p)(p^{3}-p^{2})$. This group is terminal.

\bigskip\noindent\textbf{Group 7.6.3}%

\[
\langle a,b\rangle\times_{\lbrack b,a]=[d,c]=[f,e]}\langle c,d\rangle
\times_{\lbrack b,a]=[d,c]=[f,e]}\langle e,f\rangle
\]

The number of conjugacy classes is $p^{6}+p-1$, and the automorphism group has
order $(p^{7}-p)(p^{7}-p^{6})(p^{5}-p)p^{4}(p^{3}-p)p^{2}$. This group is terminal.

\section{Order $p^{8}$}

\setcounter{subsection}{3}

\subsection{Four generator groups}

\noindent\textbf{Group 8.4.1}%

\[
\langle a,b,c,d\,|\,[b,a],\,[c,a]\rangle
\]

The number of conjugacy classes is $2p^{5}+p^{4}-2p^{3}$, and the automorphism
group has order $(p-1)^{2}(p^{2}-1)(p^{2}-p)p^{21}$. The nucleus has dimension
12. There are $9p+58+3\gcd(p-1,3)$ orbits of subspaces of the nucleus of
codimension 1, and a total of $9p+71+3\gcd(p-1,3)$ descendants of order
$p^{9}$ with exponent $p$.

\bigskip\noindent\textbf{Group 8.4.2}%

\[
\langle a,b,c,d\,|\,[b,a],\,[d,c]\rangle
\]

The number of conjugacy classes is $p^{5}+3p^{4}-2p^{3}-2p^{2}+p$, and the
automorphism group has order $2(p^{2}-1)^{2}(p^{2}-p)^{2}p^{16}$. The nucleus
has dimension 12. The number of orbits of subspaces of the nucleus of
codimension 1 is%
\begin{align*}
&  \frac{1}{2}(p^{4}+2p^{3}+7p^{2}+19p+2\gcd(p-1,3)+2\gcd(p-1,4)\\
&  +\frac{1}{2}\gcd(p-1,5)+\frac{1}{2}\gcd(p+1,5)+76).
\end{align*}
The total number of descendants of order $p^{9}$ with exponent $p$ is%
\begin{align*}
&  \frac{1}{2}(p^{4}+2p^{3}+7p^{2}+19p+2\gcd(p-1,3)+2\gcd(p-1,4)\\
&  +\frac{1}{2}\gcd(p-1,5)+\frac{1}{2}\gcd(p+1,5)+88).
\end{align*}

\bigskip\noindent\textbf{Group 8.4.3}%

\[
\langle a,b,c,d\,|\,[b,a],\,[d,b][c,a]\rangle
\]

The number of conjugacy classes is $p^{5}+2p^{4}-p^{3}-p^{2}$, and the
automorphism group has order $(p-1)(p^{2}-1)(p^{2}-p)p^{20}$. The nucleus has
dimension 12. When $p=5$ there are 310 orbits of subspace of the nucleus of
codimension 1, and when $p>5$ the number of orbits is%
\[
p^{3}+3p^{2}+12p+(p+4)\gcd(p-1,3)+2\gcd(p-1,4)+\gcd(p-1,5)+31.
\]
The total number of descendants of order $p^{9}$ with exponent $p$ is 7 more
than the number of these orbits.

\bigskip\noindent\textbf{Group 8.4.4}%

\[
\langle a,b,c,d\,|\,[d,b][c,a],\,[d,c][b,a]^{\omega}\rangle
\]

The number of conjugacy classes is $p^{5}+p^{4}-p$, and the automorphism group
has order $2(p^{4}-1)(p^{4}-p^{2})p^{16}$. The nucleus has dimension 12. The
number of orbits of subspaces of the nucleus of codimension 1 is%
\[
\frac{1}{2}(p^{4}+3p^{2}+5p+11)+\frac{1}{4}(\gcd(p^{2}+1,5)-1).
\]
The number of descendants of order $p^{9}$ with exponent $p$ is equal to the
number of these orbits.

\subsection{Five generator groups}

\noindent\textbf{Group 8.5.1}%

\[
\langle
a,b,c,d,e\,|\,[e,a],\,[c,b],\,[d,b],\,[e,b],\,[d,c],\,[e,c],\,[e,d]\rangle
\]

The number of conjugacy classes is $p^{6}+p^{5}-p^{3}$, and the automorphism
group has order $(p-1)^{2}(p^{3}-1)(p^{3}-p)(p^{3}-p^{2})p^{22}$. The nucleus
has dimension 9, and there are 8 orbits of subspaces of the nucleus of
codimension 1. There are 28 descendants of order $p^{9}$ with exponent $p$.

\bigskip\noindent\textbf{Group 8.5.2}%

\[
\langle
a,b,c,d,e\,|\,[d,a],\,[e,a],\,[d,b],\,[e,b],\,[d,c],\,[e,c],\,[e,d]\rangle
\]

The number of conjugacy classes is $p^{6}+p^{5}-p^{3}$, and the automorphism
group has order $(p^{2}-1)(p^{2}-p)(p^{3}-1)(p^{3}-p)(p^{3}-p^{2})p^{21}$. The
nucleus has dimension 8 and there are $p+4+\gcd(p-1,3)$ orbits of subspaces of
the nucleus of codimension 1. There are $2p+14+2\gcd(p-1,3)$ descendants of
order $p^{9}$ with exponent $p$.

\bigskip\noindent\textbf{Group 8.5.3}%

\[
\langle
a,b,c,d,e\,|\,[d,a],\,[e,a],\,[c,b],\,[e,b],\,[d,c],\,[e,c],\,[e,d]\rangle
\]

The number of conjugacy classes is $3p^{5}-3p^{3}+p^{2}$, and the automorphism
group has order $2(p-1)^{5}p^{23}$. The nucleus has dimension 8, and there are
$p+26$ orbits of subspaces of the nucleus of codimension 1. There are $2p+101$
descendants of order $p^{9}$ with exponent $p$.

\bigskip\noindent\textbf{Group 8.5.4}%

\[
\langle
a,b,c,d,e\,|\,[d,a]=[c,b],\,[e,a],\,[d,b],\,[e,b],\,[d,c],\,[e,c],\,[e,d]\rangle
\]

The number of conjugacy classes is $2p^{5}+p^{4}-2p^{3}$, and the automorphism
group has order $(p-1)^{2}(p^{2}-1)(p^{2}-p)p^{24}$. The nucleus has dimension
8, and there are $p+10$ orbits of subspaces of the nucleus of codimension 1.
There are $2p+37$ descendants of order $p^{9}$ with exponent $p$.

\bigskip\noindent\textbf{Group 8.5.5}%

\[
\langle
a,b,c,d,e\,|\,[c,a],\,[d,a]=[c,b],\,[e,a],\,[d,b],\,[e,b],\,[e,c],\,[e,d]\rangle
\]

The number of conjugacy classes is $2p^{5}+p^{4}-2p^{3}$, and the automorphism
group has order $(p^{2}-1)^{2}(p^{2}-p)^{2}p^{19}$. The nucleus has dimension
8, and there are $2p+10+\gcd(p-1,3)$ orbits of subspaces of the nucleus of
codimension 1. There are $5p+35+3\gcd(p-1,3)$ descendants of order $p^{9}$
with exponent $p$.

\bigskip\noindent\textbf{Group 8.5.6}%

\[
\langle a,b,c,d,e\,|\,[d,b]=[c,a]^{\omega}%
,\,[d,c]=[b,a],\,[e,a],\,[c,b],\,[e,b],\,[e,c],\,[e,d]\rangle
\]

The number of conjugacy classes is $p^{5}+2p^{4}-p^{3}-p^{2}$, and the
automorphism group has order $2(p-1)(p^{2}-1)^{2}p^{23}$. The nucleus has
dimension 8, and there are $p+5$ orbits of subspaces of the nucleus of
codimension 1. There are $2p+17$ descendants of order $p^{9}$ with exponent
$p$.

\bigskip\noindent\textbf{Group 8.5.7}%

\[
\langle a,b,c,d,e\,|\,[e,b]=[c,a][d,b]^{m}%
,[d,c]=[b,a],[e,c]=[d,b],[d,a],[e,a],[c,b],[e,d]\rangle,
\]
where $m$ is chosen so that $x^{3}+mx-1$ is irreducible over GF$(p)$.
Different choices of $m$ give isomorphic groups.

The number of conjugacy classes is $p^{5}+p^{4}-p$, and the automorphism group
has order $3(p-1)(p^{3}-1)p^{18}$. The nucleus has dimension 6, and there are
2 orbits of subspaces of the nucleus of codimension 1. The number of
descendants of order $p^{9}$ with exponent $p$ is%
\[
\frac{p^{2}+p+10+2\gcd(p-1,3)}{3}.
\]

\bigskip\noindent\textbf{Group 8.5.8}%

\[
\langle a,b,c,d,e\,|\,[d,a][c,b]^{\omega}%
,\,[e,a],\,[e,b]=[c,a],\,[d,b],\,[d,c]=[b,a],\,[e,c],\,[e,d]\rangle
\]

The number of conjugacy classes is $p^{5}+p^{4}-p$, and the automorphism group
has order $2(p-1)(p^{2}-1)p^{18}$. The nucleus has dimension 5, and there are
$p+4$ orbits of subspaces of the nucleus of codimension 1. There are $\frac
{5}{2}(p+1)+9$ descendants of order $p^{9}$ with exponent $p$.

\bigskip\noindent\textbf{Group 8.5.9}%

\[
\langle
a,b,c,d,e\,|\,[d,a],\,[e,a],\,[e,b]=[c,a],\,[d,b],\,[d,c]=[b,a],\,[e,c],\,[e,d]=[c,b]\rangle
\]

The number of conjugacy classes is $p^{5}+p^{4}-p$, and the automorphism group
has order $(p+1)(p-1)^{2}p^{16}$.The nucleus has dimension 5, and there are
$2p+6+\gcd(p-1,3)$ orbits of subspaces of the nucleus of codimension 1. There
are $4p+17+2\gcd(p-1,3)$ descendants of order $p^{9}$ with exponent $p$.

\bigskip\noindent\textbf{Group 8.5.10}%

\[
\langle
a,b,c,d,e\,|\,[d,a],\,[e,a],\,[d,b],\,[e,b]=[c,a],\,[d,c]=[b,a],\,[e,c],\,[e,d]\rangle
\]

The number of conjugacy classes is $p^{5}+2p^{4}-p^{3}-p^{2}$, and the
automorphism group has order $(p-1)(p^{2}-1)(p^{2}-p)p^{21}$. The nucleus has
dimension 6, and there are $p+5$ orbits of subspaces of the nucleus of
codimension 1. There are $2p+20$ descendants of order $p^{9}$ with exponent
$p$.

\bigskip\noindent\textbf{Group 8.5.11}%

\[
\langle
a,b,c,d,e\,|\,[d,a],\,[e,a],\,[c,b],\,[e,b]=[c,a],\,[d,c]=[b,a],\,[e,c],\,[e,d]\rangle
\]

The number of conjugacy classes is $p^{5}+2p^{4}-p^{3}-p^{2}$, and the
automorphism group has order $(p-1)^{3}p^{21}$. The nucleus has dimension 6,
and there are 9 orbits of subspaces of the nucleus of codimension 1. There are
37 descendants of order $p^{9}$ with exponent $p$.

\bigskip\noindent\textbf{Group 8.5.12}%

\[
\langle
a,b,c,d,e\,|\,[e,a],\,[c,b],\,[d,b],\,[e,c],\,[e,d],\,[d,c]=[b,a],\,[e,b]=[c,a]\rangle
\]

The number of conjugacy classes is $p^{5}+2p^{4}-p^{3}-p^{2}$, and the
automorphism group has order $(p-1)^{3}p^{19}$. The nucleus has dimension 5,
and there are 9 orbits of subspaces of the nucleus of codimension 1. There are
$p+33+\gcd(p-1,3)$ descendants of order $p^{9}$ with exponent $p$.

\bigskip\noindent\textbf{Group 8.5.13}%

\[
\langle a,b,c,d,e\,|\,[d,a],\,[e,a],\,[c,b],\,[e,b][d,b]^{\omega
}=[c,a],\,[d,c]=[b,a],\,[e,c],\,[e,d]\rangle
\]

The number of conjugacy classes is $p^{5}+2p^{4}-p^{3}-p^{2}$, and the
automorphism group has order $2(p-1)^{2}(p^{2}-1)p^{18}$. The nucleus has
dimension 6, and there are 11 orbits of subspaces of the nucleus of
codimension 1. The number of descendants of order $p^{9}$ with exponent $p$ is%
\[
\frac{p^{2}+4p+69}{2}.
\]

\bigskip\noindent\textbf{Group 8.5.14}%

\[
\langle
a,b,c,d,e\,|\,[c,b],\,[d,b],\,[e,c],\,[e,d],\,[d,c]=[b,a],\,[e,b]=[c,a],\,[e,a]=[d,a]^{\omega
}\rangle
\]

The number of conjugacy classes is $p^{5}+2p^{4}-p^{3}-p^{2}$, and the
automorphism group has order $2(p-1)(p^{2}-1)p^{17}$. The nucleus has
dimension 5, and there are 6 orbits of subspaces of the nucleus of codimension
1. There are $\frac{1}{2}(p+37)$ descendants of order $p^{9}$ with exponent
$p$.

\bigskip\noindent\textbf{Group 8.5.15}%

\[
\langle
a,b,c,d,e\,|\,[d,a],\,[e,a],\,[c,b],\,[d,b],\,[e,c],\,[d,c]=[b,a],\,[e,b]=[c,a]\rangle
\]

The number of conjugacy classes is $p^{5}+2p^{4}-p^{3}-p^{2}$, and the
automorphism group has order $2(p-1)^{3}p^{17}$.The nucleus has dimension 5,
and there are 12 orbits of subspaces of the nucleus of codimension 1. There
are $\frac{1}{2}(p+81)$ descendants of order $p^{9}$ with exponent $p$.

\bigskip\noindent\textbf{Group 8.5.16}%

\[
\langle
a,b,c,d,e\,|\,[d,a],\,[e,a],\,[c,b],\,[d,b],\,[e,b]=[c,a],\,[e,c],\,[e,d]\rangle
\]

The number of conjugacy classes is $p^{5}+3p^{4}-2p^{3}-2p^{2}+p$, and the
automorphism group has order $(p-1)^{4}p^{19}$. The nucleus has dimension 6,
and there are 18 orbits of subspaces of the nucleus of codimension 1. There
are $p+83+\gcd(p-1,3)$ descendants of order $p^{9}$ with exponent $p$.

\bigskip\noindent\textbf{Group 8.5.17}%

\[
\langle a,b,c,d,e\,|\,[e,b]=[c,a],\,[d,c]=[b,a][d,a]^{-1}%
,\,[e,a],\,[c,b],\,[d,b],\,[e,c],\,[e,d]\rangle
\]

The number of conjugacy classes is $p^{5}+3p^{4}-2p^{3}-2p^{2}+p$, and the
automorphism group has order $2(p-1)^{3}p^{18}$. The nucleus has dimension 5,
and there are $p+8$ orbits of subspaces of the nucleus of codimension 1. There
are $\frac{1}{2}(5p+63)$ descendants of order $p^{9}$ with exponent $p$.

\bigskip\noindent\textbf{Group 8.5.18}%

\[
\langle
a,b,c,d,e\,|\,[e,c][e,b],\,[c,a],\,[d,a],\,[e,a],\,[c,b],\,[d,b],\,[e,d]\rangle
\]

The number of conjugacy classes is $p^{5}+4p^{4}-3p^{3}-3p^{2}+2p$, and the
automorphism group has order $6(p-1)^{4}p^{18}$. The nucleus has dimension 6,
and there are $12$ orbits of subspaces of the nucleus of codimension 1. The
number of descendants of order $p^{9}$ with exponent $p$ is%
\[
\frac{1}{6}(p^{2}+10p+241+2\gcd(p-1,3)).
\]

\bigskip\noindent\textbf{Group 8.5.19}%

\[
\langle
a,b,c,d,e\,|\,[d,a],\,[e,a],\,[c,b],\,[d,b],\,[e,b],\,[d,c],\,[e,c]\rangle
\]

The number of conjugacy classes is $2p^{5}+2p^{4}-3p^{3}-p^{2}+p$, and the
automorphism group has order $(p-1)(p^{2}-1)^{2}(p^{2}-p)^{2}p^{17}$. The
nucleus has dimension 7, and there are $11$ orbits of subspaces of the nucleus
of codimension 1. There are $40$ descendants of order $p^{9}$ with exponent
$p$.

\bigskip\noindent\textbf{Group 8.5.20}%

\[
\langle
a,b,c,d,e\,|\,[d,c]=[b,a],\,[c,a],\,[d,a],\,[e,a],\,[d,b],\,[e,b],\,[e,d]\rangle
\]

The number of conjugacy classes is $2p^{5}+p^{4}-2p^{3}$, and the automorphism
group has order $(p-1)^{4}p^{23}$. The nucleus has dimension 7, and there are
$15$ orbits of subspaces of the nucleus of codimension 1. There are $67$
descendants of order $p^{9}$ with exponent $p$.

\bigskip\noindent\textbf{Group 8.5.21}%

\[
\langle
a,b,c,d,e\,|\,[d,c]=[b,a],\,[c,a],\,[d,a],\,[e,a],\,[c,b],\,[d,b],\,[e,d]\rangle
\]

The number of conjugacy classes is $2p^{5}+p^{4}-2p^{3}$, and the automorphism
group has order $(p-1)^{2}(p^{2}-1)(p^{2}-p)p^{20}$. The nucleus has dimension
6, and there are $10$ orbits of subspaces of the nucleus of codimension 1.
There are $45$ descendants of order $p^{9}$ with exponent $p$.

\bigskip\noindent\textbf{Group 8.5.22}%

\[
\langle
a,b,c,d,e\,|\,[d,c]=[b,a],\,[c,a],\,[d,a],\,[e,a],\,[c,b],\,[d,b],\,[e,b]\rangle
\]

The number of conjugacy classes is $2p^{5}+p^{4}-2p^{3}$, and the automorphism
group has order $(p^{2}-1)^{2}(p^{2}-p)^{2}p^{17}$. The nucleus has dimension
5, and there are $5$ orbits of subspaces of the nucleus of codimension 1.
There are $p+17$ descendants of order $p^{9}$ with exponent $p $.

\subsection{Six generator groups}

For all these groups we take the generators to be $a,b,c,d,e,f$, and we just
give the relations, with the class two and exponent $p$ conditions understood.

\bigskip\noindent\textbf{Group 8.6.1}%

\begin{align*}
&  \lbrack c,b],\,[d,a],\,[d,b],\,[d,c],\,[e,a],\,[e,b],\,[e,c],\\
&  \lbrack e,d],\,[f,a],\,[f,b],\,[f,c],\,[f,d],\,[f,e]
\end{align*}

The number of conjugacy classes is $2p^{6}-p^{4}$ and the order of the
automorphism group is $(p-1)(p^{2}-1)(p^{2}-p)(p^{3}-1)(p^{3}-p)(p^{3}%
-p^{2})p^{23}$. The nucleus has dimension 5, and there are $5$ orbits of
subspaces of the nucleus of codimension 1. There are $28$ descendants of order
$p^{9}$ with exponent $p$.

\bigskip\noindent\textbf{Group 8.6.2}%

\begin{align*}
&  [c,b],\,[d,a],\,[d,b]=[b,a],\,[d,c],\,[e,a],\,[e,b],\,[e,c],\\
&  [e,d],\,[f,a],\,[f,b],\,[f,c],\,[f,d],\,[f,e]
\end{align*}

The number of conjugacy classes is $p^{6}+2p^{5}-p^{4}-2p^{3}+p^{2}$ and the
order of the automorphism group is $2(p^{2}-1)^{3}(p^{2}-p)^{3}p^{20}$. The
nucleus has dimension 4, and there are $2$ orbits of subspaces of the nucleus
of codimension 1. There are $19$ descendants of order $p^{9}$ with exponent
$p$.

\bigskip\noindent\textbf{Group 8.6.3}%

\begin{align*}
&  [c,b],\,[d,a],\,[d,b]=[c,a],\,[d,c],\,[e,a],\,[e,b],\,[e,c],\\
&  [e,d],\,[f,a],\,[f,b],\,[f,c],\,[f,d],\,[f,e]
\end{align*}

The number of conjugacy classes is $p^{6}+p^{5}-p^{3}$ and the order of the
automorphism group is $(p-1)(p^{2}-1)^{2}(p^{2}-p)^{2}p^{24}$. The nucleus has
dimension 4, and there are $2$ orbits of subspaces of the nucleus of
codimension 1. There are $21$ descendants of order $p^{9}$ with exponent $p$.

\bigskip\noindent\textbf{Group 8.6.4}%

\begin{align*}
&  [c,b],\,[d,a],\,[d,b]=[c,a],\,[d,c]=[b,a]^{\omega}%
,\,[e,a],\,[e,b],\,[e,c],\\
&  [e,d],\,[f,a],\,[f,b],\,[f,c],\,[f,d],\,[f,e]
\end{align*}

The number of conjugacy classes is $p^{6}+p^{4}-p^{2}$ and the order of the
automorphism group is $2(p^{4}-1)(p^{4}-p^{2})(p^{2}-1)(p^{2}-p)p^{20}$. The
nucleus has dimension 4, and there is only 1 orbit of subspaces of the nucleus
of codimension 1. There are $7$ descendants of order $p^{9}$ with exponent $p$.

\bigskip\noindent\textbf{Group 8.6.5}%

\begin{align*}
&  [c,b],\,[d,a],\,[d,b],\,[d,c],\,[e,a],\,[e,b],\,[e,c],\\
&  [e,d]=[b,a],\,[f,a],\,[f,b],\,[f,c],\,[f,d],\,[f,e]
\end{align*}

The number of conjugacy classes is $p^{6}+p^{5}-p^{3}$ and the order of the
automorphism group is $(p-1)^{3}(p^{2}-1)(p^{2}-p)p^{22}$. The nucleus has
dimension 2, and there are $2$ orbits of subspaces of the nucleus of
codimension 1. There are $31$ descendants of order $p^{9}$ with exponent $p$.

\bigskip\noindent\textbf{Group 8.6.6}%

\begin{align*}
&  [c,b],\,[d,a],\,[d,b]=[c,a],\,[d,c],\,[e,a],\,[e,b],\,[e,c],\\
&  [e,d]=[b,a],\,[f,a],\,[f,b],\,[f,c],\,[f,d],\,[f,e]
\end{align*}

The number of conjugacy classes is $p^{6}+p^{4}-p^{2}$ and the order of the
automorphism group is $(p-1)^{2}(p^{2}-1)(p^{2}-p)p^{21}$. The nucleus has
dimension 1, and there are $10$ descendants of order $p^{9}$ with exponent $p
$.

\bigskip\noindent\textbf{Group 8.6.7}%

\begin{align*}
&  [c,a],\,[c,b],\,[d,a],\,[d,b],\,[e,a],\,[e,b],\,[e,c],\\
&  [e,d],\,[f,a],\,[f,b],\,[f,c],\,[f,d],\,[f,e]=[b,a]
\end{align*}

The number of conjugacy classes is $p^{6}+p^{5}-p^{4}+p^{3}-2p+1$ and the
order of the automorphism group is $(p^{4}-1)(p^{4}-p^{3})(p^{2}-1)^{2}%
(p^{2}-p)p^{13}$. The nucleus has dimension 2, and there is only one orbit of
subspaces of the nucleus of codimension 1. There are $10$ descendants of order
$p^{9}$ with exponent $p$.

\bigskip\noindent\textbf{Group 8.6.8}%

\begin{align*}
&  [c,a],\,[c,b],\,[d,a],\,[d,b],\,[e,a],\,[e,b],\,[e,c],\\
&  [e,d],\,[f,a],\,[f,b],\,[f,c],\,[f,d],\,[f,e]=[b,a][d,c]
\end{align*}

The number of conjugacy classes is $p^{6}+3p^{3}-2p^{2}-3p+2$ and the order of
the automorphism group is $6(p^{2}-1)^{3}(p^{2}-p)p^{14}$. This group is terminal.

\bigskip\noindent\textbf{Group 8.6.9}%

\begin{align*}
&  [c,b],\,[d,a],\,[d,b]=[c,a],\,[d,c],\,[e,a],\,[e,b],\,[e,c],\\
&  [e,d],\,[f,a],\,[f,b],\,[f,c],\,[f,d],\,[f,e]=[b,a]
\end{align*}

The number of conjugacy classes is $p^{6}+2p^{3}-p^{2}-2p+1$ and the order of
the automorphism group is $(p^{2}-1)^{2}(p^{2}-p)^{2}p^{15}$. This group is terminal.

\bigskip\noindent\textbf{Group 8.6.10}%

\begin{align*}
&  [c,b],\,[d,a],\,[d,b]=[c,a],\,[d,c],\,[e,a],\,[e,b],\,[e,c],\\
&  [e,d],\,[f,a],\,[f,b],\,[f,c],\,[f,d],\,[f,e]=[c,a]
\end{align*}

The number of conjugacy classes is $p^{6}+p^{5}-p^{4}+p^{2}-p$ and the order
of the automorphism group is $(p^{2}-1)^{2}(p^{2}-p)^{2}p^{20}$. The nucleus
has dimension 2, and there is only 1 orbit of subspaces of the nucleus of
codimension 1. There are $13$ descendants of order $p^{9}$ with exponent $p$.

\bigskip\noindent\textbf{Group 8.6.11}%

\begin{align*}
&  [c,b],\,[d,a],\,[d,b]=[c,a],\,[d,c]=[b,a]^{\omega}%
,\,[e,a],\,[e,b],\,[e,c],\\
&  [e,d],\,[f,a],\,[f,b],\,[f,c],\,[f,d],\,[f,e]=[b,a]
\end{align*}

The number of conjugacy classes is $p^{6}+p^{3}-p$ and the order of the
automorphism group is $2(p^{4}-1)(p^{3}-p^{2})(p^{2}-1)p^{13}$. This group is terminal.

\bigskip\noindent\textbf{Group 8.6.12}%

\begin{align*}
&  [c,b],\,[d,a],\,[d,b],\,[d,c],\,[e,a],\,[e,b],\,[e,c],\\
&  [e,d]=[b,a],\,[f,a],\,[f,b],\,[f,c],\,[f,d]=[c,a],\,[f,e]
\end{align*}

The number of conjugacy classes is $p^{6}+p^{4}-p^{2}$ and the order of the
automorphism group is $(p^{2}-1)^{2}(p^{2}-p)^{2}p^{18}$. This group is terminal.

\bigskip\noindent\textbf{Group 8.6.13}%

\begin{align*}
&  [b,a],\,[d,a],\,[e,a][c,a],\,[f,a],\,[c,b],\,[d,b]=[c,a],\,[e,b],\\
&  [f,b][c,a]^{2},\,[d,c],\,[e,c],\,[e,d]=[f,c],\,[f,d],\,[f,e]=[c,a][f,c]
\end{align*}

The number of conjugacy classes is $p^{6}+p^{3}-p$ and the order of the
automorphism group is $(p-1)(p^{2}-1)(p^{2}-p)p^{19}$. This group is terminal.

\bigskip\noindent\textbf{Group 8.6.14}%

\begin{align*}
&  [c,b],\,[d,a],\,[d,b]=[c,a],\,[d,c],\,[e,a],\,[e,b],\,[e,c],\\
&  [e,d]=[b,a],\,[f,a],\,[f,b]=[c,a]^{m},\,[f,c]=[b,a],\,[f,d],\,[f,e]=[c,a],
\end{align*}
where $x^{3}-mx+1$ is irreducible over GF$(p)$. (Different choices of $m$ give
isomorphic groups.)

The number of conjugacy classes is $p^{6}+p^{2}-1$ and the order of the
automorphism group is $3(p^{6}-1)(p-1)p^{15}$. This group is terminal.

\subsection{Seven generator groups}

\noindent\textbf{Group 8.7.1}%

\[
\langle a,b\rangle\times\langle c\rangle\times\langle d\rangle\times\langle
e\rangle\times\langle f\rangle\times\langle g\rangle.
\]

The number of conjugacy classes is $p^{7}+p^{6}-p^{5}$, and the automorphism
group has order $(p^{8}-p^{6})(p^{8}-p^{7})(p^{6}-p)(p^{6}-p^{2})(p^{6}%
-p^{3})(p^{6}-p^{4})(p^{6}-p^{5}) $. The nucleus has dimension 2, and there is
only one orbit of subspaces of the nucleus of dimension 1. There are $6$
descendants of order $p^{9}$ with exponent $p$.

\bigskip\noindent\textbf{Group 8.7.2}%

\[
\langle a,b\rangle\times_{\lbrack b,a]=[d,c]}\langle c,d\rangle\times\langle
e\rangle\times\langle f\rangle\times\langle g\rangle.
\]

The number of conjugacy classes is $p^{7}+p^{4}-p^{3}$, and the automorphism
group has order $(p^{8}-p^{4})(p^{8}-p^{7})(p^{6}-p^{4})p^{5}(p^{4}%
-p)(p^{4}-p^{2})(p^{4}-p^{3}) $. This group is terminal.

\bigskip\noindent\textbf{Group 8.7.3}%

\[
\langle a,b\rangle\times_{\lbrack b,a]=[d,c]=[f,e]}\langle c,d\rangle
\times_{\lbrack b,a]=[d,c]=[f,e]}\langle e,f\rangle\times\langle g\rangle.
\]

The number of conjugacy classes is $p^{7}+p^{2}-p$, and the automorphism group
has order $(p^{8}-p^{2})(p^{8}-p^{7})(p^{6}-p^{2})p^{5}(p^{4}-p^{2}%
)p^{3}(p^{2}-p)$. This group is terminal.

\end{document}